\documentclass[letterpaper,11pt]{amsart}

\usepackage{amssymb}
\usepackage{eucal}
\usepackage{array}
\usepackage{epsfig}
\usepackage[all]{xy}
\CompileMatrices \numberwithin{equation}{section}
\SelectTips{cm}{}

\title{More about vanishing cycles and mutation}
\author{Paul Seidel}
\date{October 5, 2000.}

\renewcommand{\subsection}{\vspace{1ex}
\hspace{-\parindent}%
\stepcounter{subsection}%
({\sc \thesection\alph{subsection}})}

\newcommand{\includefigure}[3]{%
\begin{figure}[#3]
\begin{center}
\epsfig{file=#2} \\ \caption{\label{fig:#1}}
\end{center}
\end{figure}}

\newenvironment{myitemize}%
{\begin{itemize}\itemsep0.5em} {\end{itemize}}

\newcommand{\R}{\mathbb{R}}
\newcommand{\Z}{\mathbb{Z}}
\newcommand{\C}{\mathbb{C}}
\newcommand{\Q}{\mathbb{Q}}

\newcommand{\iso}{\cong}           
\newcommand{\smooth}{C^\infty}
\newcommand{\CP}[1]{\C {\mathrm P}^{#1}}
\newcommand{\RP}[1]{\R {\mathrm P}^{#1}}

\newcommand{\Rgeq}{\R^{\scriptscriptstyle \geq 0}}
\newcommand{\Rleq}{\R^{\scriptscriptstyle \leq 0}}
\newcommand{\suchthat}{\; | \;}

\newcommand{\id}{\mathrm{id}}

\newcommand{\im}{\mathrm{im}}
\renewcommand{\ker}{\mathrm{ker}}

\newcommand{\Hom}{\mathrm{Hom}}

\renewcommand{\o}{\omega}
\renewcommand{\O}{\Omega}

\newcommand{\Symp}{\mathrm{Symp}}
\newcommand{\Diff}{\mathrm{Diff}}

\theoremstyle{plain}
\newtheorem{thm}{Theorem}[section]
\newtheorem{theorem}[thm]{Theorem}

\newtheorem{lemma}[thm]{Lemma}
\newtheorem{prop}[thm]{Proposition}

\newtheorem{conjecture}[thm]{Conjecture}
\newtheorem{question}[thm]{Question}
\newtheorem{defn}[thm]{Definition}

\newtheorem{remark}[thm]{Remark}
\newtheorem{remarks}[thm]{Remarks}
\newtheorem{example}[thm]{Example}

\theoremstyle{remark}
\newtheorem*{additional}{Additional references}
\newtheorem*{acknowledgments}{Acknowledgments}

\newcommand{\Sympe}{\Symp^e}
\newcommand{\crit}{\mathrm{crit}}
\newcommand{\Phirel}{\Phi_{rel}}
\newcommand{\A}{\mathcal A}
\newcommand{\B}{\mathcal B}
\newcommand{\Cone}{\mathrm{Cone}}
\newcommand{\tL}{\widetilde{L}}
\newcommand{\tphi}{\tilde{\phi}}
\newcommand{\ttau}{\tilde{\tau}}

\newcommand{\Lag}{{\mathit{Lag}}^\rightarrow}
\renewcommand{\hom}{\mathit{hom}}
\newcommand{\Tw}{\mathrm{Tw}}
\newcommand{\Ob}{\mathrm{Ob}}
\newcommand{\JJ}{\mathcal{J}}
\newcommand{\RR}{\mathcal{R}}
\renewcommand{\SS}{\mathcal{S}}
\newcommand{\sing}{\mathrm{sing}}

\newcommand{\Morse}{\mathit{Morse}}
\newcommand{\fund}[1]{[\hspace{-0.15em}[ #1 ]\hspace{-0.15em}]}

\begin{document}
\maketitle

\section{Introduction}
Not surprisingly, this is a sequel to ``Vanishing cycles and mutation''
\cite{seidel00}. Notions from that paper will be used freely later on, but for
starters a short recapitulation seems appropriate. Let $D \subset \C$ be the
closed unit disc. An {exact Morse fibration} $\pi: E^{2n+2} \rightarrow D$ can
be described roughly as a family of $2n$-dimensional exact symplectic manifolds
parametrized by $D$, with singularities modelled on Morse-type critical points
of holomorphic functions. Such objects are suitable for developing
Picard-Lefschetz theory in a symplectic context. Suppose that $(E,\pi)$ comes
equipped with a relative Maslov map $\delta_{E/D}$, which is basically a
trivialization of $(\Lambda^{n+1}_{\C} TE)^{\otimes 2}$. Then one can associate
to it an algebraic invariant, the derived directed Fukaya category
$D^b\Lag(\widetilde{\Gamma})$, which encodes the symplectic geometry of the
vanishing cycles of $(E,\pi)$. It is constructed in three steps:

\begin{myitemize}
\item[(i)]
Fix a base point $z_0 \in \partial D$, for definiteness $z_0 = -i$; the fibre
$M = E_{z_0}$ is an exact symplectic manifold. Then make an admissible choice
of paths \cite[Figure 2]{seidel00} from $z_0$ to the critical values of $\pi$.
This determines a {distinguished basis of vanishing cycles} in $M$, which is an
ordered family $\Gamma = (L_1,\dots,L_m)$ of embedded exact Lagrangian spheres
$L_i \subset M$ (there is one more bit of information, so-called framings of
the $L_i$, which we ignore for the moment). We will sometimes call such
families Lagrangian configurations in $M$. Any other distinguished basis,
arising from a different choice of paths, can be obtained from $\Gamma$ through
a sequence of Hurwitz moves.

\item[(ii)]
$\delta_{E/D}$ induces a Maslov map $\delta_M$ on $M$, thus allowing one to
speak of graded Lagrangian submanifolds in $M$. By choosing gradings in an
arbitrary way, one lifts $\Gamma$ to a graded Lagrangian configuration
$\widetilde{\Gamma} = (\tL_1,\dots,\tL_m)$. Next one introduces the {directed
Fukaya category} $\A = \Lag(\widetilde{\Gamma})$, which is an
$A_\infty$-category with objects $\tL_1,\dots,\tL_m$. The spaces of morphisms
in $\A$ are essentially given by the natural cochain complexes underlying the
Floer cohomology groups $HF^*(\tL_i,\tL_k)$, but the ordering of the objects
also plays a role in the definition, which is what the word ``directed'' refers
to.

\item[(iii)]
The {derived category} $D^b(\A)$ of $\A$ is no longer an $A_\infty$-category
but rather a category in the ordinary sense; in fact it is a triangulated
category, linear over $\Z/2$ since we are using $\Z/2$-coefficients. Its
objects are twisted complexes in $\A$, which are $A_\infty$-analogues of chain
complexes, and morphisms are chain homotopy classes of maps in a suitable
sense. Unlike $\A$ itself, $D^b(\A)$ is independent of the choice of $\Gamma$
(and of $\widetilde{\Gamma}$), hence is an invariant of $(E,\pi)$ and
$\delta_{E/D}$. The reason for this is a relation between Hurwitz moves of
$\Gamma$ and the purely algebraic notion of mutations of $\A$.
\end{myitemize}
All this is explained in more detail, although still without proofs, in
\cite{seidel00}. The present paper adds computations and some more far-reaching
conjectures. It is divided into four more or less independent parts. The first
part, formed by Sections \ref{sec:dimension-zero}--\ref{sec:examples},
introduces some classes of examples. Section \ref{sec:dimension-zero} concerns
the toy model case when the fibre dimension is zero. Section \ref{sec:examples}
explains how to construct exact Morse fibrations from algebro-geometric
objects, namely Lefschetz pencils or isolated hypersurface singularities. Also
included is a discussion of the mirror manifold of $\CP{2}$.

The subject of the second part, spanning Sections
\ref{sec:tqft2}--\ref{sec:donaldson}, is Hochschild cohomology. The basic
observation is that like $D^b(\A)$, the Hochschild cohomology $HH^*(\A,\A)$ is
an invariant of $(E,\pi)$ and $\delta_{E/D}$. One would expect this to have a
more geometric meaning. Following a suggestion of Donaldson we propose a
conjecture in this direction, which relates Hochschild cohomology to the global
monodromy of the fibration. As a preliminary to this, Section \ref{sec:tqft2}
enlarges the topological quantum field theory framework from \cite[Section
3]{seidel00} to include the Floer cohomology of symplectic automorphisms.

The third part, Sections \ref{sec:morse}--\ref{sec:real}, concerns exact Morse
fibrations with a real involution. In that situation, and under certain
additional assumptions, the category $D^b(\A)$ should have a description in
terms of Morse theory on the real part. We state a precise conjecture, which is
a theorem at least in the lowest nontrivial dimension $n = 1$.

Finally, Section \ref{sec:induction} is a reflection on the aspect of
dimensional reduction inherent in the theory. By this we mean that $\A$, which
is defined in terms of pseudo-holomorphic curves in the $2n$-dimensional
manifold $M$, produces invariants of the $(2n+2)$-dimensional object $(E,\pi)$.
In an attempt to push this further, we give a conjectural algorithm for
determining the Floer cohomology of certain Lagrangian spheres in $E$ through
computations in $D^b(\A)$. This, if correct, might be the first true
application of homological algebra methods in symplectic geometry.

\begin{acknowledgments}
This paper owes much to Donaldson; several of the ideas presented here arose
during conversations with him. Givental and Kontsevich provided helpful
explanations concerning mirror symmetry.
\end{acknowledgments}

\newpage

\section{Branched covers\label{sec:dimension-zero}}
\subsection{}
In \cite{seidel00} it was tacitly assumed that all symplectic manifolds were of
dimension $>0$ (and that assumption will be resumed as soon as this section is
finished). The zero-dimensional case is rather trivial, but it still requires
some minor adaptations, which we now explain. Let's start with a basic
dictionary of zero-dimensional symplectic geometry:

\begin{center}
\begin{tabular}{l||l}
 compact symplectic manifold & finite set $M$ \\ \hline
 Lagrangian sphere & subset $L \subset M$, $|L| = 2$
 \\ \hline
 Floer cohomology & $HF(L_1,L_2) = (\Z/2)^{L_1 \cap L_2}$
 \\ \hline
 Dehn twist & \begin{minipage}{16em} \vspace{2pt}
 \raggedright
 $\tau_L: M \rightarrow M$ swapping the two
 points of $L$, leaving the rest fixed \vspace{2pt}
 \end{minipage}
 \\ \hline
 exact Morse fibration & generic branched cover $\pi: E \rightarrow S$
\end{tabular}
\end{center}
Here a {\bf generic branched cover} is a proper map $\pi: E \rightarrow S$
between oriented surfaces with boundary, with $\pi^{-1}(\partial S) = \partial
E$, having finitely many critical points, which are locally modelled on $\pi(z)
= z^2$; and no two critical points may lie in the same fibre.

The TQFT framework deserves a brief mention, if only for its simplicity. Let
$S$ be a compact oriented surface with boundary, $\Sigma \subset \partial S$ a
finite set of marked points, and $(E,\pi)$ a generic branched cover over $S^* =
S \setminus \Sigma$. Suppose moreover that we have a submanifold $Q \subset
\pi^{-1}(\partial S^*)$ such that $\pi|Q: Q \rightarrow \partial S^*$ is a
double cover; this corresponds to what in \cite{seidel00} was called a
Lagrangian boundary condition for $(E,\pi)$. The behaviour of $Q$ near a point
$\zeta \in \Sigma$ can be described by a pair $(L_{\zeta,+},L_{\zeta,-})$ of
Lagrangian zero-spheres in some fibre of $E$ close to $\zeta$. A continuous
section $u: S^* \rightarrow E$ of $\pi$ which satisfies $u(\partial S^*)
\subset Q$ singles out a point of $L_{\zeta,+} \cap L_{\zeta,-}$ for each
$\zeta$, and hence determines an element of $\bigotimes_\zeta HF(L_{\zeta,+},
L_{\zeta,-})$. Taking the sum over all sections defines an invariant
\begin{equation} \label{eq:rel}
\Phirel(E,\pi,Q) \in \bigotimes_{\zeta \in \Sigma} HF(L_{\zeta,+},L_{\zeta,-}).
\end{equation}
For trivial reasons, these invariants satisfy the gluing law formulated in
\cite[Section 3]{seidel00}. One can use them to find the correct
zero-dimensional analogues of various familiar maps in Floer theory. For
instance, the pair-of-pants product $HF(L_2,L_3) \otimes HF(L_1,L_2)
\rightarrow HF(L_1,L_3)$ is derived from an invariant \eqref{eq:rel} with $S =
D$, $\Sigma = \{\text{\em three points}\}$ and $E = S^* \times M$, together
with a boundary condition given by $L_1,L_2,L_3$ over the three connected
components of $\partial S^*$. Working this out explicitly shows that the
product takes a basis element $\bar{x} \otimes \bar{y} \in HF(L_2,L_3) \otimes
HF(L_1,L_2)$, $x \in L_2 \cap L_3$ and $y \in L_1 \cap L_2$, to $\bar{x} \in
HF(L_1,L_3)$ if $x = y$, and to zero otherwise. There is the same long exact
sequence in Floer cohomology as in \cite[Theorem 3.3]{seidel00} and again, the
more fundamental fact underlying it is a quasi-isomorphism
\[
\Cone(a: HF(L,L_2) \otimes HF(L_1,L) \rightarrow HF(L_1,L_2))
\xrightarrow{(h,b)} HF(L_1,\tau_L(L_2)).
\]
However, while $a$ and $b$ are as in higher dimensions, $h$ is the product
$HF(L,L_2) \otimes HF(L_1,L) \iso HF(L,\tau_L(L_2)) \otimes HF(L_1,L)
\rightarrow HF(L_1,\tau_L(L_2))$. Of course, both the result and proof are now
elementary!

The next topic are gradings. We summarize their zero-dimensional theory in
another dictionary:

\begin{center}
\begin{tabular}{l||l}
 \begin{minipage}{9em} \vspace{2pt}
 Maslov map on a \newline symplectic manifold \vspace{2pt}
 \end{minipage}
 & $\delta_M: M \rightarrow S^1$
 \\ \hline
 \begin{minipage}{9em} \vspace{2pt} \raggedright
 graded Lagrangian sphere \vspace{2pt}
 \end{minipage}
 &
 \begin{minipage}{18em} \vspace{2pt} \raggedright
 $L \subset M$ with a map $\tL: L \rightarrow \R$ such
 that $\exp(2\pi i \tL) = \delta_M|L$ \vspace{2pt}
 \end{minipage}
 \\ \hline
 \begin{minipage}{9em} \vspace{2pt} \raggedright
 graded Floer cohomology \vspace{2pt}
 \end{minipage}
 &
 \begin{minipage}{18em} \vspace{2pt} \raggedright
 in $HF^*(\tL_1,\tL_2)$, the degree of $\bar{x}$ is
 $\tL_2(x) - \tL_1(x)$
 \vspace{2pt}
 \end{minipage}
 \\ \hline
 shifting the grading &
 $\tL[\sigma] = \tL - \sigma$, $\sigma \in \Z$
 \\ \hline
 \begin{minipage}{9em} \vspace{2pt} \raggedright
 graded symplectic automorphism \vspace{2pt}
 \end{minipage}
 &
 \begin{minipage}{18em} \vspace{2pt}
 $\phi: M \rightarrow M$ with a map $\tphi: M \rightarrow \R$ such that
 $\exp(2\pi i \phi(x)) = \delta_M(x)/\delta_M(\phi^{-1}(x))$
 \vspace{2pt}
 \end{minipage}
 \\ \hline
 relative Maslov map &
 \begin{minipage}{18em} \vspace{2pt} \raggedright
 on $\pi: E \rightarrow S$, a trivialization
 $\delta_{E/S}$ of $TE^{\otimes 2} \otimes (\pi^*TS)^{\otimes -2}$
 \end{minipage}
\end{tabular}
\end{center}
The action of a graded symplectic automorphism on graded Lagrangian
submanifolds is defined by $(\tphi \tL)(x) = \tL(\phi^{-1}(x)) + \tphi(x)$. A
grading $\tL$ determines a preferred grading $\ttau_{\tL}$ of the associated
Dehn twist, characterized by
\[
\ttau_{\tL}(\tL) = \tL[1] \quad \text{ and } \quad \ttau_{\tL}(x) = 0 \text{
for } x \notin L;
\]
compare \cite[Equation (4)]{seidel00} for the first property. Graded Lagrangian
configurations $\widetilde{\Gamma} = (\tL_1,\dots,\tL_m)$ are defined as one
would expect. Their Hurwitz moves are
\begin{itemize}
\item
$\widetilde{\Gamma} \rightsquigarrow (\tL_1[\sigma_1],\dots,\tL_m[\sigma_m])$
for $\sigma_1,\dots,\sigma_m \in \Z$;
\item
$\widetilde{\Gamma} \rightsquigarrow c\widetilde{\Gamma} =
(\ttau_{\tL_1}(\tL_2),\dots,\ttau_{\tL_1}(\tL_m),\tL_1)$;
\item
$\widetilde{\Gamma} \rightsquigarrow r\widetilde{\Gamma} =
(\tL_1,\dots,\tL_{m-2},\ttau_{\tL_{m-1}}(\tL_m),\tL_{m-1})$.
\end{itemize}
Let $(E,\pi)$ be a generic branched cover over $D$, with a relative Maslov map
$\delta_{E/D}$. A drastically simplified version of the usual Picard-Lefschetz
argument produces from it a graded (with respect to the induced Maslov map
$\delta_M$) Lagrangian configuration in $M = E_{z_0}$, unique up to Hurwitz
equivalence.

\begin{remark}
One can observe here two minor differences with respect to the
positive-dimensional situation. Firstly, $\ttau_{\tL}$ depends on $\tL$;
secondly, the Hurwitz equivalence class of a graded configuration contains
information that is lost if one forgets the gradings. Both phenomena go back to
the obvious fact that there is a $\Z^2$ ambiguity in the choice of grading for
a Lagrangian zero-sphere, as opposed to $\Z$ in higher dimensions.
\end{remark}

The directed Fukaya category $\A = \Lag(\widetilde{\Gamma})$ attached to a
zero-dimensi\-onal graded configuration has objects $\tL_i$. The morphisms
$hom_\A(\tL_i,\tL_k)$ are $HF^*(\tL_i,\tL_k)$ for $i<k$, $\Z/2 \cdot
\id_{\tL_i}$ for $i = k$, and zero for $i>k$. The nontrivial compositions
$HF^*(\tL_j,\tL_k) \otimes HF^*(\tL_i,\tL_j) \rightarrow HF^*(\tL_i,\tL_k)$,
$i<j<k$, are given by the pair-of-pants product. Even though $\A$ is an
ordinary $\Z/2$-linear and $\Z$-graded category, we prefer to regard it as a
directed $A_\infty$-category where all composition maps $\mu_\A^d$ of order $d
\neq 2$ vanish. The details of the proof are slightly different than in the
higher-dimensional case, but the main result of \cite{seidel00} remains true,
which is that Hurwitz moves of $\widetilde{\Gamma}$ give rise to mutations of
$\A$; and one arrives at the familiar conclusion that
$D^b\Lag(\widetilde{\Gamma})$ is an invariant of $(E,\pi)$, $\delta_{E/D}$.

Before going on to concrete examples, we need to recall some algebraic
terminology. Let $T$ be a quiver (an oriented graph). The path category
${\mathcal P}T$ has one object for each vertex, and the space of morphisms
between two objects is the $\Z/2$-vector space freely generated by all paths in
$T$ going from one vertex to the other. One can further divide the morphism
spaces by some two-sided ideal, in order to kill certain designated morphisms.
This is usually referred to as describing a category by a ``quiver with
relations''. If the quiver is directed, which means that its vertices are
numbered $1,\dots,m$ such that there are no paths from the $i$-th one to the
$k$-th one unless $i<k$, ${\mathcal P}T$ and its quotients can again be seen as
special cases of directed $A_\infty$-categories, with morphisms only in degree
zero and vanishing composition maps of order $\neq 2$. Note that in this
situation our notion of derived category (defined through twisted complexes)
agrees with the classical one (defined by considering the categories as
algebras, and taking chain complexes of right modules over them).

\subsection{}
Take the following quiver with $m$ vertices, with the arrows oriented in an
arbitrary way:
\begin{equation} \label{eq:bgp}
\xymatrix{
 {\bullet} \ar[r] &
 {\bullet} & \ar[l]
 {\bullet} \ar[r] &
 {\dots} &
 {\bullet} \ar[l]
}
\end{equation}
and let $\A$ be its path category. It is a classical result that $D^b(\A)$ is
independent of the orientation of the arrows; what we will do is to explain
this geometrically. First of all, one can find a configuration of Lagrangian
zero-spheres $\Gamma = (L_1,\dots,L_m)$ in $M = \{1,\dots,m+1\}$ which, when
lifted to a graded configuration in the trivial way ($\delta_M = 0$ and $\tL_i
= 0$), has $\A$ as its directed Fukaya category. The way to do that is best
explained by an example:
\begin{equation} \label{eq:makebasis}
\begin{split}
& \xymatrix@R=-0.25em{
 \bullet \ar[r] &
 \bullet \ar[r] &
 \bullet &
 \bullet \ar[l] &
 \bullet \ar[l] &
 \bullet \ar[l] \ar[r] &
 \bullet
 \\
 {\hspace{-1em} \{1,2\} \hspace{-1em}} &
 {\hspace{-1em} \{1,3\} \hspace{-1em}} &
 {\hspace{-1em} \{1,4\} \hspace{-1em}} &
 {\hspace{-1em} \{5,4\} \hspace{-1em}} &
 {\hspace{-1em} \{6,4\} \hspace{-1em}} &
 {\hspace{-1em} \{7,4\} \hspace{-1em}} &
 {\hspace{-1em} \{7,8\} \hspace{-1em}}
}
\\[0.5em]
\Gamma & = (\{1,2\},\{1,3\},\{7,4\},\{6,4\},\{5,4\},\{1,4\},\{7,8\}).
\end{split}
\end{equation}
Next one constructs an $(m+1)$-fold generic branched cover $\pi: E \rightarrow
D$ such that $\Gamma$ is one of its distinguished bases of vanishing cycles.
Inspection of \eqref{eq:makebasis} and an Euler characteristic computation show
that the total space $E$ is connected and in fact a disc. It is a consequence
of the classification of branched covers that there is just one such $(E,\pi)$
up to isomorphism, and moreover it admits precisely one homotopy class of
relative Maslov maps. Therefore all categories $D^b(\A)$ arise from the same
geometric situation, hence must be equivalent on grounds of the general theory.

\subsection{}
Let $\pi: E \rightarrow D$ be a double cover branched along $2g+1 \geq 5$
points, so that the total space is a genus $g \geq 2$ surface with one boundary
component. Write $\iota$ for the nontrivial covering transformation, which is a
hyperelliptic involution of $E$. There is a unique homotopy class of relative
Maslov maps which are invariant under $\iota$; take $\delta_{E/D}$ in that
class. All vanishing cycles of $(E,\pi)$ are the same, and if we choose their
gradings in the most obvious way, the resulting directed Fukaya category $\A_g$
can be described by the quiver with relations
\[
\xymatrix{
 {\bullet} \ar@/^1pc/[r]^{a_1} \ar@/_1pc/[r]_{b_1} &
 {\bullet} \ar@/^1pc/[r]!<-1em,0em>^{a_2} \ar@/_1pc/[r]!<-1em,0em>_{b_2} &
 {\bullet\; \dots\; \bullet}
 \ar@/^1pc/[]!<1em,0em>;[r]^{a_{2g}}
 \ar@/_1pc/[]!<1em,0em>;[r]_{b_{2g}} &
 {\bullet}
} \qquad \qquad
\begin{cases}
 b_{i+1}a_i = 0, & \\
 a_{i+1}b_i = 0.
\end{cases}
\]
We will now make some further remarks concerning $\A_g$ and its derived
category. These may seem a bit unmotivated, but they will be at least partially
put into context later on; see Remark \ref{th:zero-spheres}. The following
definition is taken from \cite{seidel-thomas99}:

\begin{defn} \label{def:spherical}
Let ${\mathcal C}$ be a triangulated category, linear over $\Z/2$ and such that
the spaces $\Hom^*_{\mathcal C}(-,-)$ are finite-dimensional. $C \in
\Ob\,{\mathcal C}$ is called spherical of dimension $(n+1)$ if
$\Hom^*_{\mathcal C}(C,C) \iso H^*(S^{n+1};\Z/2)$ and the composition
\[
\Hom^*_{\mathcal C}(X,C) \otimes \Hom^{n+1-*}_{\mathcal C}(C,X) \rightarrow
\Hom^{n+1}_{\mathcal C}(C,C) \iso \Z/2
\]
is a nondegenerate pairing for any $X \in \Ob\,{\mathcal C}$.
\end{defn}

Under some additional assumptions on ${\mathcal C}$, which are satisfied for
derived categories of directed $A_\infty$-categories, one can associate to any
spherical object $C$ an exact self-equivalence of ${\mathcal C}$, the twist
functor $T_C$, which is well-defined up to isomorphism. If $C_1,C_2$ are
spherical with $\Hom^*_{\mathcal C}(C_1,C_2) = 0$ then their twist functors
commute; and if $\Hom^*_{\mathcal C}(C_1,C_2)$ is one-dimensional, one gets a
braid relation $T_{C_1}T_{C_2}T_{C_1} \iso T_{C_2}T_{C_1}T_{C_2}$.

Denote the objects of $\A_g$ by $(X^1,\dots,X^{2g+1})$. The twisted complexes
$C_i = (X^i[1] \oplus X^{i+1}, \left(\begin{smallmatrix} 0 & 0 \\ a_i+b_i & 0
\end{smallmatrix}\right) )$, $1 \leq i \leq 2g$, are spherical objects of
dimension one in $D^b(\A_g)$, and moreover
\[
 \dim\, \Hom^*_{D^b(\A_g)}(C_i,C_k) =
 \begin{cases}
 0 & |i-k| \geq 2, \\
 1 & |i-k| = 1.
 \end{cases}
\]
This implies that $T_{C_1},\dots,T_{C_{2g}}$ generate an action (in a weak
sense) of the braid group $B_{2g+1}$ on $D^b(\A_g)$. Now let $D^{per}(\A_g)$ be
the category defined like $D^b(\A_g)$ but using $\Z/2$-graded twisted
complexes. This is triangulated, and its two-fold shift functor $[2]$ is
isomorphic to the identity; we call it the periodic derived category of $\A_g$.
One can adapt the theory of spherical objects and twist functors to such
categories. In our case, this means that there are $C_i \in D^{per}(\A_g)$
defined as before, with the same properties. In addition there is now another
spherical object
\[
 C_0 = (X^1 \oplus X^2[1] \oplus X^3 \oplus X^4[1],
 \mbox{\smaller $\begin{pmatrix}
 0         & 0   & 0   & 0 \\
 a_1       & 0   & 0   & 0 \\
 0         & b_2 & 0   & 0 \\
 b_3b_2b_1 & 0   & a_3 & 0
 \end{pmatrix}$} )
\]
The dimension of $\Hom^*_{D^{per}(\A_g)}(C_0,C_i)$ is one if $i = 4$ and zero
for all other $i>0$. The diagram representing this situation,
\[
\xymatrix{
  & & & {C_0} \ar@{-}[d] \\
  {C_1} \ar@{-}[r] &
  {C_2} \ar@{-}[r] &
  {C_3} \ar@{-}[r] &
  {C_4} \ar@{-}[r] &
  {C_5} \ar@{-}[r] &
  \dots \ar@{-}[r] &
  {C_{2g}}
}
\]
has the same structure as the configuration of curves in Wajnryb's presentation
of $R_{g,1}$, the mapping class group of a genus $g$ surface with one marked
point \cite{wajnryb83}; and one can ask

\begin{question}
Let ${\mathcal C} \subset D^{per}(\A_g)$ be the full triangulated subcategory
generated by the $C_i$. Do $T_{C_0},\dots,T_{C_{2g}}$ generate an action of
$R_{g,1}$ on ${\mathcal C}$?
\end{question}

\begin{additional}
The classical paper about the quiver \eqref{eq:bgp} is by Bernstein, Gelfand,
and Ponomarev \cite{bernstein-gelfand-ponomarev72}. The category $D^b(\A_g)$ is
closly related to those studied in \cite{khovanov-seidel98}, even though the
quiver presentation looks different. Those papers do not consider periodic
derived categories and therefore miss out on the additional object $C_0$.
\end{additional}

\section{Examples from algebraic geometry\label{sec:examples}}
\subsection{}
Let $X$ be a smooth projective variety, $\xi \rightarrow X$ an ample line
bundle, and $\sigma_0,\sigma_1$ two holomorphic sections of $\xi$ which
generate a Lefschetz pencil of hypersurfaces $Y_z = \{ x \in X \suchthat
\sigma_0(x)/\sigma_1(x) = z\}$, $z \in \CP{1} = \C \cup \{\infty\}$. Suppose,
in addition to the Lefschetz condition, that $Y_\infty$ is smooth. What we want
to look at is, in principle, the holomorphic Morse function
\[
\sigma_0/\sigma_1 : X \setminus Y_\infty \longrightarrow \C.
\]
To obtain additional symplectic data, one chooses a metric on $\xi$, with
corresponding connection $A$, such that $\o_X = (i/2\pi)F_A$ is a K{\"a}hler
form. By a standard construction, $A$ and $\sigma_1$ give rise to a one-form
$\theta_X$ defined on $X \setminus Y_\infty$, such that $d\theta_X = \o_X$.
This still doesn't fit the definition of an exact Morse fibration, making some
further modifications necessary.

Let $Z = Y_0 \cap Y_\infty$ be the base locus of the pencil, $G = \{ (z,x) \in
\CP{1} \times X \suchthat \sigma_0(x)/\sigma_1(x) = z\}$ its graph, and $p: G
\rightarrow \CP{1}$ the projection. $(G,p)$ is a Lefschetz fibration with
fibres $Y_z$, and contains a trivial subfibration $\CP{1} \times Z \subset G$.
Let $\O_G \in \Omega^2(G)$ be the pullback of $\o_X$; it is closed,
nondegenerate on the vertical tangent spaces $\ker(Dp) \subset TG$, and hence
defines a symplectic connection away from the critical points of $p$. Moreover,
the restriction of this connection to $\CP{1} \times Z \subset G$ is trivial.
Similarly, pulling back $\theta_X$ gives a one-form $\Theta_G$ defined on the
complement of $(\{\infty\} \times Y_\infty) \cup (\CP{1} \times Z) \subset G$,
with $d\Theta_G = \Omega_G$.

After possibly rescaling $\sigma_0$, we may assume that $Y_z$ is smooth for
$|z| \geq 1$. Let $U \subset Y_0$ be a small open neighbourhood of $Z$. Using
the symplectic parallel transport in radial directions one constructs an
embedding $\Xi$,
\[
\xymatrix{
 {D \times U} \ar[rr]^-{\Xi} \ar[dr] &&
 {p^{-1}(D) \subset G} \ar[dl]^p \\
 & {D} &
}
\]
with $\Xi\, |\, D \times Z$ and $\Xi\, |\,\{0\} \times U$ the obvious
inclusions, which has the property that $\Xi^*\Omega_G - \o_X|U$ vanishes on
each diameter $[-z;z] \times U$, $z \in S^1$. Let $h \in \smooth_c(U,[0;1])$ be
a function which is $\equiv 1$ on some smaller neighbourhood $U' \subset U$ of
$Z$. Take the map $H: p^{-1}(D) \rightarrow p^{-1}(D)$ which is the identity
outside $\im(\Xi)$ and satisfies $H\,\Xi(z,x) = \Xi((1-h(x))z,x)$.

Set $E = p^{-1}(D) \setminus \Xi(D \times U'')$ for an even smaller open
neighbourhood $U'' \subset U'$ of $Z$, which we choose such that $\partial U''$
is smooth. Together with $\pi = p|E$, $\Omega = (H^*\Omega_G|E) \in
\Omega^2(E)$, $\Theta = (H^*\Theta_G|E) \in \Omega^1(E)$, and the given complex
structure near the critical points of $\pi$, this forms an exact Morse
fibration over $D$. Indeed, the required trivialization near $\partial_hE$ is
given by $\Xi$ itself; and the fact that $H$ contracts along the same radial
directions which are used to define $\Xi$ ensures that $\Omega$ equals
$\Omega_G$ on $\ker(D\pi) \subset TE$, hence is nondegenerate there.

Suppose now that the canonical bundle $K_X$ satisfies
\begin{equation} \label{eq:relative-maslov}
K_X^{-2} \iso \xi^a \, \text{ holomorphically, for some $a \in \Z$.}
\end{equation}
Then $\sigma_1^a$ gives a trivialization of $K_X^{-2}$ over $X \setminus
Y_\infty$, determined up to multiplication with a constant. Since the
projection $G \rightarrow X$ identifies $E$ with a compact subset of $X
\setminus Y_\infty$, and $\O$ agrees with the pullback of $\o_X$ except in a
small neighbourhood of $\partial_hE$, one gets, at least up to homotopy, a
preferred relative Maslov map $\delta_{E/D}$. The category
$D^b\Lag(\widetilde{\Gamma})$ associated to $(E,\pi)$ and $\delta_{E/D}$
actually depends only on $X$ and on the Lefschetz pencil, and not on the other
choices made during the construction of $E$. That is because different choices
lead to exact Morse fibrations which are deformation equivalent in a suitable
sense. For the same reason, if $\xi$ is such that any two generic sections
define a Lefschetz pencil, the category is an invariant of the pair $(X,\xi)$.
\includefigure{degree2}{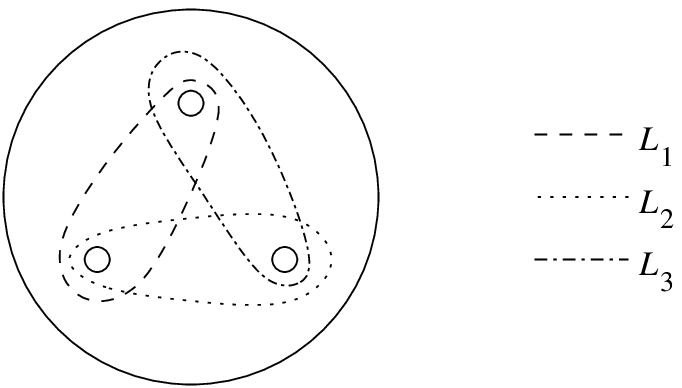}{hb}%

As a concrete example take $X = \CP{2}$, $\xi = {\mathcal O}(2)$. The fibre $M
= E_{z_0}$ of $(E,\pi)$ is, topologically, $\CP{1}$ with four small discs
removed. The relative Maslov map $\delta_{E/D}$ coming from
\eqref{eq:relative-maslov} induces a Maslov map $\delta_M$ on $M$, which is
characterized by having the same behaviour over each component of $\partial M$.
Figure \ref{fig:degree2} shows a distinguished basis $\Gamma = (L_1,L_2,L_3)$
of vanishing cycles. One can read off directly that the directed Fukaya
category $\A = \Lag(\widetilde{\Gamma})$, for a suitable choice of gradings, is
described by the quiver with relations
\begin{equation} \label{eq:cp2-quiver}
\xymatrix{
 {\bullet} \ar@/^1pc/[rr]^{a_1} \ar@/_1pc/[rr]_{b_1} &&
 {\bullet} \ar@/^1pc/[rr]^{a_2} \ar@/_1pc/[rr]_{b_2} &&
 {\bullet}
} \qquad
\begin{cases}
 a_2a_1 = b_2b_1, & \\
 b_2a_1 = a_2b_1. &
\end{cases}
\end{equation}
Namely, the intersections $L_i \cap L_k$, $i<k$, consist of two points which
are essential, in the sense that they cannot be removed by an isotopy in $M$.
This implies that $hom_\A(\tL_i,\tL_k)$ is two-dimensional and $\mu^1_\A = 0$.
And there are four triangles in $M$ whose sides (in positive order) map to
$(L_1,L_2,L_3)$, giving rise to the nonzero products $\mu^2_\A(a_2,a_1)$,
$\mu^2_\A(b_2,a_1)$, $\mu_\A^2(a_2,b_1)$, $\mu^2_\A(b_2,b_1)$.

\begin{remark}
Suppose that \eqref{eq:relative-maslov} holds but with $a \in \Q$, which is to
say $K_X^{-2c} \iso \xi^b$ for coprime $b,c$. Then $\sigma_1^b$ defines a
$c$-sheeted cyclic covering of $X \setminus Y_\infty$, hence also a covering of
$E$. The total space of the latter can be made into an exact Morse fibration
which then has a canonical relative Maslov map, so that our theory can be
applied to it. The $\Z/c$-action on $E$ by covering transformations induces a
self-equivalence $\sigma$ of $D^b\Lag(\widetilde{\Gamma})$ which satisfies
$\sigma^c \iso [1]$. One can view this as a category with ``fractional
gradings'', thinking of $\Hom(X,\sigma^iY)$ as ``the group $\Hom^{i/c}(X,Y)$ of
morphisms of degree $i/c$''.
\end{remark}

\subsection{}
We need to extend the previous discussion to a degenerate situation. Let $X
\subset \CP{N}$ be a projective variety with singular set $X^\sing$, and
$\sigma_0,\sigma_1 \in (\C^{N+1})^*$ two linear forms which generate a pencil
$(Y_z)$ of hyperplane sections of $X$; so the line bundle we are looking at is
$\xi = {\mathcal O}_X(1)$. Suppose that the following conditions are satisfied:
$X^\sing \subset Y_\infty \setminus Y_0$; the base locus $Z = Y_0 \cap
Y_\infty$ is smooth; and $(Y_z)_{z \neq \infty}$ satisfies the same
nondegeneracy conditions as a Lefschetz pencil. Then the associated graph $p: G
\rightarrow \CP{1}$ is well-behaved except over $\infty \in \CP{1}$. The
construction made above, in which the fibre at infinity plays no role, goes
through without any changes, yielding an exact Morse fibration $(E,\pi)$ over
$D$. If we suppose moreover that $X$ is locally a complete intersection, so
that there is a well-defined canonical bundle $K_X$, then condition
\eqref{eq:relative-maslov} again ensures the existence of preferred relative
Maslov maps $\delta_{E/D}$.
\includefigure{mirrorp2}{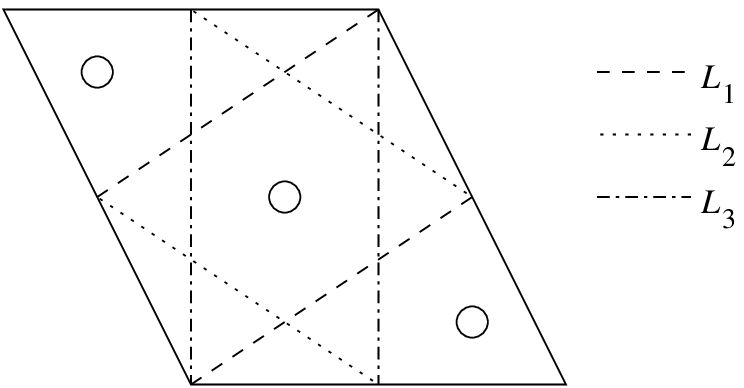}{hb}%

The example which we have in mind is $X = \{ x \in \CP{3} \suchthat x_0^3 =
x_1x_2x_3\}$, $\sigma_0(x) = x_1+x_2+x_3$, $\sigma_1(x) = x_0$. One can
identify $X \setminus Y_\infty$ with the affine hypersurface $x_1x_2x_3 = 1$,
and then $\sigma_0/\sigma_1$ becomes the function $x_1 + x_2 + x_3$, which
appears in Givental's work \cite[Theorem 5]{givental94} as the object mirror
dual to the projective plane. The next result exhibits the mirror phenomenon
within our framework.

\begin{prop} \label{th:mirror}
Let $(E,\pi)$ be the exact Morse fibration constructed from
$(X,\sigma_0,\sigma_1)$, with relative Maslov map $\delta_{E/D}$ obtained from
\eqref{eq:relative-maslov}. Then for any graded distinguished basis of
vanishing cycles ones has
\[
D^b\Lag(\widetilde{\Gamma}) \iso D^b Coh({\mathbb P}^2)
\]
where $D^b Coh({\mathbb P}^2)$ is the derived category, in the classical sense,
of coherent sheaves on the projective plane over the field $\Z/2$.
\end{prop}

The fibre $M = E_{z_0}$ is a torus with three small discs removed. Figure
\ref{fig:mirrorp2} shows a particular distinguished basis of vanishing cycles
$\Gamma = (L_1,L_2,L_3)$. In that picture, the Maslov map $\delta_M$ is that
given by the flat structure of $M$. One can therefore choose gradings such that
the groups $CF^*(\tL_i,\tL_k)$, $i<k$, are concentrated in degree zero. Then
$\A = \Lag(\widetilde{\Gamma})$ is given by the quiver with relations
\[
\xymatrix{
 {\bullet} \ar@/^1.5pc/[rr]^{a_1} \ar[rr]^{a_2} \ar@/_1.5pc/[rr]_{a_3} &&
 {\bullet} \ar@/^1.5pc/[rr]^{b_1} \ar[rr]^{b_2} \ar@/_1.5pc/[rr]_{b_3} &&
 {\bullet}
} \qquad
\begin{cases}
 b_i a_j = b_j a_i & i \neq j, \\
 b_i a_i = 0.
\end{cases}
\]
A theorem of Beilinson \cite{beilinson78} says that $D^b(\A) \iso D^b
Coh({\mathbb P}^2)$. Since $D^b(\A)$ remains the same for all other
distinguished bases, the Proposition is proved. The same approach can be used
for the mirror partner of $\CP{1} \times \CP{1}$, with analogous results; I
have not checked any further cases.

\begin{remarks}
(i) It is maybe helpful to mention one way of making drawings like Figure
\ref{fig:mirrorp2}. By projecting to one coordinate, each fibre
$(\sigma_0/\sigma_1)^{-1}(z)$ of $\sigma_0/\sigma_1: X \setminus Y_\infty
\rightarrow \C$ can be represented as a double covering of $\C^*$ branched over
three points. The vanishing cycles are preimages of paths in $\C^*$ joining two
branch points, and one finds them by looking at how the branch points come
together as one varies $z$.

(ii) One sees from Figure \ref{fig:mirrorp2} that the Dehn twists $\tau_{L_i}$
can be represented by affine maps, at least away from a small neighbourhood of
the missing discs. It follows that any Lagrangian configuration $\Gamma' =
(L_1',L_2',L_3')$ which is Hurwitz equivalent to $\Gamma$ consists of curves
$L_i'$ isotopic in $M$ to straight lines. An easy argument using the topology
of $E$ 
shows that no two $L_i'$ can ever be isotopic to each other. Therefore
$HF^*(\tL_i',\tL_k')$, $i<k$, is concentrated in a single degree, for any
choice of gradings. This is related to the fact that mutations of the standard
exceptional collection in $D^b Coh({\mathbb P}^2)$ remain strongly exceptional;
see \cite{bondal-polishchuk94} for much more about this subject.
\end{remarks}

\subsection{}
Let $f \in \C[x_1,\dots,x_{n+1}]$ be a polynomial with an isolated critical
point at the origin. A Morsification of $f$ is a smooth family of polynomials
$(f_t)_{0 \leq t < \rho}$ with $f_0 = f$ and such that the critical points of
$f_t$, $t>0$, which lie near the origin are nondegenerate. For sufficiently
small $0<\epsilon$, $0<\delta \ll \epsilon$, and $0 < t \ll \delta$, one finds
that
\[
X = \{ x \in \C^{n+1} \suchthat |x| \leq \epsilon, \; |f_t(x)| \leq \delta\}
\]
is a manifold with corners, and $f_t$ a holomorphic Morse function on it. One
can equip $X$ with the standard symplectic form on $\C^{n+1}$, and its
canonical bundle has a standard trivialization obtained from the constant
holomorphic volume form. After some modifications, this becomes an exact Morse
fibration $(E,\pi)$ over $D$ together with a canonical homotopy class of
relative Maslov maps $\delta_{E/D}$. We will not explain the details since the
procedure is similar to that for Lefschetz pencils. A deformation equivalence
argument shows that the resulting category $D^b\Lag(\widetilde{\Gamma})$ is
independent of all choices, and even of the Morsification, so that it is an
invariant of $f$.

\begin{prop} \label{th:ade}
Let $f \in \C[x_1,x_2]$ be an ADE singularity. Then the category
$D^b\Lag(\widetilde{\Gamma})$ associated to it is equivalent to the derived
path category of the Dynkin quiver of the same type, with an arbitrary
orientation of the arrows.
\end{prop}

A'Campo \cite[pp.\ 13--17]{acampo75b} has constructed particularly nice
Morsifications and distinguished bases $\Gamma = (L_1,\dots,L_m)$ for these
singularities. These have the properties that (i) there is a bijection between
the $L_i$ and the vertices of the corresponding Dynkin diagram, such that $L_i
\cap L_k$ consists of a single point if the $i$-th and $k$-th vertex are
connected by a line, and is empty otherwise; (ii) there are no $i<k<l$ with
$L_i \cap L_k \neq \emptyset$ and $L_k \cap L_l \neq \emptyset$. Figure
\ref{fig:e6} shows the $E_6$ case. One can choose gradings $\tL_i$ such that
$HF^*(\tL_i,\tL_k)$ is concentrated in degree zero for all $i<k$. Then
$\Lag(\widetilde{\Gamma})$ is the path category of the Dynkin quiver for a
particular orientation, which is such that each vertex is either a sink or a
source. Changing the orientation does not affect the derived category of the
path category, which completes the proof.
\includefigure{e6}{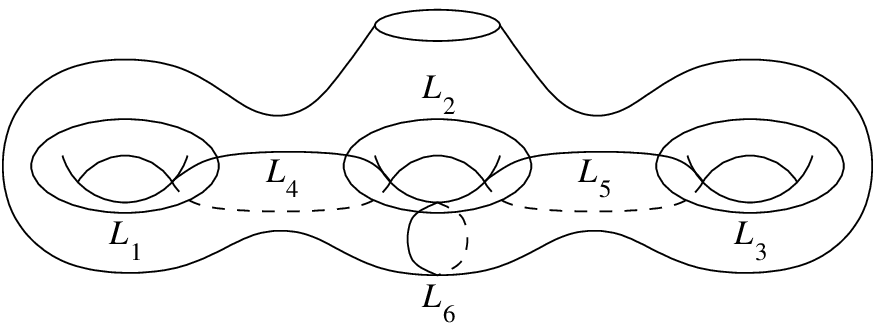}{ht}%

We mention two more expected properties of the categories
$D^b\Lag(\widetilde{\Gamma})$ as invariants of singularities. Strictly speaking
these are conjectures, but the proofs should not be difficult.
\begin{myitemize}
\item ($\mu$-invariance)
They should remain the same when one deforms $f$ inside its $\mu = \text{\em
const.}$ stratum. Moreover, if $f$ is adjacent to $g$ then the category
belonging to $f$ should contain that of $g$ as a full triangulated subcategory.
\item (Stabilization)
The categories associated to $f(x)$ and $\hat{f}(x,y) = f(x) + y^2_1 + \dots +
y_q^2$ should be equivalent, for any $q$.
\end{myitemize}
For instance, the reader may have noticed that the categories of
$A_m$-singula\-ri\-ties in Proposition \ref{th:ade} are the same as those
arising from certain branched covers, which were considered in Section
\ref{sec:dimension-zero}b. This is an extreme case of the stabilization
property, since the branched covers can be viewed as Morsifications of the
one-dimensional singularities $f(x) = x^{m+1}$.

\begin{additional}
Hori, Iqbal and Vafa \cite{hori-iqbal-vafa00} discuss mirror symmetry for Fano
varieties from a physics point of view. The construction of exact Morse
fibrations from singularities is similar to the definition of symplectic
monodromy in \cite[Section 6]{khovanov-seidel98}.
\end{additional}

\section{Floer cohomology for automorphisms\label{sec:tqft2}}
\subsection{}
Let $(M,\o,\theta)$ be an exact symplectic manifold (which is, as always,
assumed to have contact type boundary). Floer cohomology associates to each
$\phi \in \Sympe(M)$ a pair of vector spaces $HF(\phi,-)$ and $HF(\phi,+)$.
These correspond to two ways of treating the fixed points near $\partial M$:
let $H \in \smooth(M,\R)$ be a function with $H|\partial M \equiv 0$, whose
Hamiltonian flow $(\phi_t^H)$ is equal to the Reeb flow on $\partial M$. To
define $HF(\phi,-)$ and $HF(\phi,+)$ one perturbs $\phi$ to $\phi \circ
\phi_t^H$ for some small $t<0$, respectively $t>0$. The difference is measured
by a long exact sequence
\begin{equation} \label{eq:boundary-term}
\xymatrix{
 {HF(\phi,-)} \ar[r] & {HF(\phi,+)} \ar[r] & {H^*(\partial M;\Z/2).}
 \ar@/^1.5em/[ll]_{ }
}
\end{equation}
For $\phi = \id_M$ this reduces to the usual cohomology exact sequence, with
$HF(\id_M,-) \iso H^*(M,\partial M;\Z/2)$ and $HF(\id_M,+) \iso H^*(M;\Z/2)$.
Other properties of $HF(\phi,\pm)$ are its invariance under isotopies within
$\Sympe(M)$; conjugation invariance, $HF(\phi_2\phi_1,\pm) \iso
HF(\phi_1\phi_2,\pm)$; and Poincar{\'e} duality, $HF(\phi,\pm) \iso
HF(\phi^{-1},\mp)^\vee$. If $M$ is equipped with a Maslov map $\delta_M$, a
grading $\tphi$ determines a $\Z$-grading on Floer cohomology, denoted by
$HF^*(\tphi,\pm)$. This satisfies $HF^*(\tphi[\sigma],\pm) \iso
HF^{*-\sigma}(\tphi,\pm)$ for $\sigma \in \Z$.

Let $S$ be a closed oriented surface together with a finite set $\Sigma \subset
S$ of marked points, and $(E,\pi) = (E,\pi,\Omega,\Theta,J_0,j_0)$ an exact
Morse fibration over $S^* = S \setminus \Sigma$, with smooth fibres isomorphic
to $M$. Suppose that around each $\zeta \in \Sigma$ we have oriented local
coordinates $\psi_\zeta: D \rightarrow S$, such that there are commutative
diagrams
\begin{equation} \label{eq:local-structure}
\xymatrix{
 {\Rleq \times T(\phi_\zeta)} \ar[d] \ar[rr]^-{\Psi_{\zeta}} &&
 {E} \ar[d]^{\pi} \\
 {D^*} \ar[rr]^-{\psi_{\zeta}|D^*} &&
 {S^*.}
}
\end{equation}
Here $T(\phi_\zeta)$ is the mapping torus of $\phi_\zeta \in \Sympe(M)$,
meaning the quotient of $\R \times M$ by $(t,x) \sim (t-1,\phi_\zeta(x))$. The
map from $\Rleq \times T(\phi_\zeta)$ to $D^* = D \setminus \{0\}$ is
$(s,[t,x]) \mapsto \exp\,2\pi(s+it)$. The pullback of $\o$ to $\R \times M$ is
invariant under the $\Z$-action, hence descends to a two-form $\o_{\phi_\zeta}$
on $T(\phi_\zeta)$; and a choice of function $K \in \smooth_c(M \setminus
\partial M,\R)$ with $dK = \phi^*\theta - \theta$ yields a one-form
$\theta_{\phi_\zeta}$ on the same space. $\Psi_{\zeta}$ is an isomorphism of
differentiable fibre bundles between $\Rleq \times T(\phi_\zeta)$ and
$E|\psi_\zeta(D^*)$, satisfying $\Psi_\zeta^*\Omega = \o_{\phi_\zeta}$,
$\Psi_\zeta^*\Theta =\theta_{\phi_\zeta}$. By and large, this means that the
behaviour of $(E,\pi)$ around $\zeta$ is described by the exact symplectic
automorphism $\phi_\zeta$. One then has a relative invariant
\begin{equation} \label{eq:relative2}
\Phirel(E,\pi) \in \bigotimes_{\zeta \in \Sigma} HF(\phi_\zeta,\pm),
\end{equation}
where the signs $\pm$ can be chosen arbitrarily, subject to the restriction
that there must be at least one $\zeta$ labeled $+$ in each connected component
of $S$; in particular $\Sigma$ may not be empty. Changing an existing $-$ to a
$+$ gives a new relative invariant which is the image of the old one under the
map $HF(\phi_\zeta,-) \rightarrow HF(\phi_\zeta,+)$ from
\eqref{eq:boundary-term}. To define the relative invariant, one chooses a
suitable complex structure on $S^*$ and almost complex structure on $E$, and
considers the moduli spaces of pseudo-holomorphic sections.

\begin{remark}
The asymmetry between $+$ and $-$ can be explained as follows. By definition of
an exact Morse fibration, one has a trivialization near $\partial E$. However,
the almost complex structures $J$ on $E$ used here do not have the obvious
product structure with respect to this trivialization; this is made necessary
by the perturbation of $\phi_\zeta$ near $\partial M$ which defines
$HF(\phi_\zeta,\pm)$. The convexity of $\partial E$ with respect to
$J$-holomorphic sections now becomes a more delicate matter, and it is there
that the sign question appears.
\end{remark}

Suppose that $\zeta \neq \zeta'$ are points in $\Sigma$ with $\phi_{\zeta'} =
(\phi_{\zeta})^{-1}$, and look at a relative invariant $\Phirel(E,\pi)$ in
which $\zeta,\zeta'$ have different signs. The surface $\overline{S}$ obtained
by connect summing together $\zeta,\zeta'$ comes with marked points
$\overline{\Sigma}$ inherited from $\Sigma \setminus \{\zeta,\zeta'\}$, and
with an exact Morse fibration over $\overline{S} \setminus \overline{\Sigma}$.
The relative invariant of this new fibration (keeping the signs as before) can
be computed, if it is well-defined, by applying the Poincar{\'e} duality
pairing $HF(\phi_\zeta,\pm) \otimes HF(\phi_{\zeta'},\mp) \rightarrow \Z/2$ to
$\Phirel(E,\pi)$.

This setup can be generalized and unified with that in \cite[Section
3]{seidel00}. In the generalization one considers oriented compact surfaces $S$
with boundary, with marked points $\Sigma \subset S$ which may lie on the
boundary or in the interior, together with an exact Morse fibration $(E,\pi)$
over $S^* = S \setminus \Sigma$ having a Lagrangian boundary condition $Q$. The
structure of $(E,\pi)$ near a marked point in the interior remains as before,
while near boundary marked points $(E,\pi,Q)$ is as in \cite{seidel00}. One
then gets relative invariants of a mixed kind,
\begin{equation} \label{eq:mixed}
 \Phirel(E,\pi,Q) \in \bigotimes_{\zeta \in \Sigma \cap \partial S}
 HF(L_{\zeta,+},L_{\zeta,-}) \otimes \bigotimes_{\zeta \in \Sigma \setminus
 \partial S} HF(\phi_\zeta,\pm).
\end{equation}
The restriction that there must be at least one $+$ now applies only to those
connected components of $S$ which have no boundary; and there are two different
gluing formulas, for boundary and interior marked points respectively.

\subsection{}
We can now explain several kinds of maps between Floer cohomology groups as
special cases of the relative invariants \eqref{eq:relative2} and
\eqref{eq:mixed}.

\begin{myitemize}
\item
Take an exact symplectic manifold $M$ and $\phi_1,\phi_2 \in \Sympe(M)$. Let $S
= S^2$ and $\Sigma = \{\text{\em 3 points}\}$. There is a unique exact Morse
fibration $(E,\pi)$ over $S^*$ which has no critical points and is even flat
(locally trivial), such that the symplectic monodromy around the three missing
points is respectively, $\phi_1$, $\phi_2$ and $(\phi_1\phi_2)^{-1}$. By
marking the points with suitable signs and using Poincar{\'e} duality, one
obtains from $\Phirel(E,\pi)$ a map $HF(\phi_1,\pm) \otimes HF(\phi_2,+)
\rightarrow HF(\phi_1\phi_2,\pm)$, which is the ``pair-of-pants'' product in
the Floer cohomology of symplectic automorphisms.

\item
Let $L \subset M$ be an exact framed Lagrangian sphere. Then there is an exact
Morse fibration $(E,\pi)$ over $D$ with a single critical point, such that the
monodromy around $\partial D$ is $\tau_L$. If one extends this to $\C = S^2
\setminus \{\infty\}$ in such a way that it is flat near $\infty$, one gets an
invariant $\Phirel(E,\pi) \in HF(\tau_L^{-1},+)$.

\item
Take the pair-of-pants product with $\phi_1 = \phi\tau_L$ and $\phi_2 =
\tau_L^{-1}$, and plug in the distinguished element of $HF(\tau_L^{-1},+)$
constructed above. This gives a canonical map $HF(\phi\tau_L,\pm) \rightarrow
HF(\phi,\pm)$, which we call $c$. Thanks to the gluing formula, $c$ can also be
described as the relative invariant of a certain exact Morse fibration over
$S^2 \setminus \{\text{\em 2 points}\}$ with a single critical point.

\item
For any $\phi \in \Sympe(M)$, $\Rleq \times T(\phi^{-1})$ is an exact Morse
fibration over $D^*$. Denote by $(E,\pi)$ its restriction to $D \setminus
\{0,i\}$. Parallel transport along $\partial D$ starting at the base point $z_0
= -i$ yields a trivialization $E\,|\,\partial D \setminus \{i\} \iso (\partial
D \setminus \{i\}) \times M$. Given an exact Lagrangian submanifold $L \subset
M$, one can introduce the Lagrangian boundary condition which in this
trivialization is $(\partial D \setminus \{i\}) \times L$. Because of the
monodromy around $0$, $\Phirel(E,\pi,Q)$ takes values in $HF(\phi^{-1},\pm)
\otimes HF(\phi(L),L)$. Equivalently, one can see it as a map $d: HF(\phi,\pm)
\rightarrow HF(\phi(L),L)$.
\end{myitemize}

\begin{theorem} \label{th:exact-sequence2}
Let $(M,\o,\theta)$ be an exact symplectic manifold, $L \subset M$ an exact
framed Lagrangian sphere, and $\phi \in \Sympe(M)$ an exact automorphism.
Suppose that $2c_1(M,L) \in H^2(M,L)$ is zero. Then there is a long exact
sequence, with maps $c$ and $d$ as defined above,
\begin{equation} \label{eq:sequence2}
\xymatrix{
 {HF(\phi\tau_L,\pm)} \ar[r]^-{c} &
 {HF(\phi,\pm)} \ar[d]^-{d} \\ &
 {HF(\phi(L),L).} \ar@/^1pc/[ul]
}
\end{equation}
\end{theorem}

If $M$ has a Maslov map $\delta_M$, and $\phi$, $L$ have gradings $\tphi$,
$\tL$ then, with respect to the natural gradings of the groups in
\eqref{eq:sequence2}, the maps $c,d$ have degree zero while the remaining one
has degree one. The proof of the exact sequence is similar to that of
\cite[Theorem 3.3]{seidel00}.

\begin{example}
By taking $\phi = \id_M$ resp.\ $\phi = \tau_L^{-1}$ one gets exact sequences
\begin{align*}
 \xymatrix{
 {HF(\tau_L,+)} \ar[r]
 & {H^*(M;\Z/2)} \ar[d] \\ &
 {H^*(L;\Z/2),} \ar@/^1pc/[ul]
 }
 \qquad
 \xymatrix{
 {H^*(M;\Z/2)} \ar[r]
 & {HF(\tau_L^{-1},+)} \ar[d] \\ &
 {H^*(L;\Z/2).} \ar@/^1pc/[ul]
}
\end{align*}
In the case when $M$ is a surface, this fits in with the results of direct
computation in \cite{seidel96b}, which are that $HF(\tau_L,+) \iso
H^*(M,L;\Z/2)$ and $HF(\tau_L^{-1},+) \iso H^*(M \setminus L;\Z/2)$. In fact,
the same holds in all dimensions, as one can show by taking a closer look at
the map $d$ in the exact sequence.
\end{example}
\includefigure{skein2}{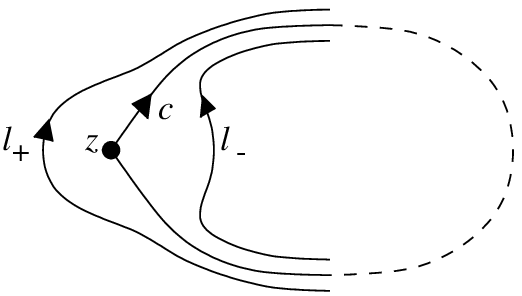}{hb}%

Now take an exact Morse fibration $(E,\pi)$ with arbitrary base $S$. To any
loop $l: S^1 \rightarrow S \setminus S^\crit$ one can associate its symplectic
monodromy $\rho_l$ and hence Floer cohomology groups, which are invariant under
homotopies of $l$ that do not pass across $S^\crit$. We write $HF(l,\pm)$
instead of $HF(\rho_l,\pm)$, in parallel with the notation $HF(c)$ adopted at
the end of \cite[Section 3]{seidel00} for paths $c: [0;1] \rightarrow S$ with
$c^{-1}(S^\crit) = \{0;1\}$. Theorem \ref{th:exact-sequence2} implies that if
one has loops $l_-,l_+$ and a path $c$ as in Figure \ref{fig:skein2}, there is
a long exact sequence
\[
\xymatrix{
 {HF(l_-)} \ar[r] & {HF(l_+)} \ar[r] & {HF(c).} \ar@/^1.5pc/[ll]^{ }
}
\]
The connection with the original formulation is made by the Picard-Lefschetz
theorem: $l_+$ differs from $l_-$ by a clockwise turn around $z \in S^\crit$,
so that $\rho_{l_+}^{-1}\rho_{l_-}$ is the Dehn twist along the vanishing cycle
arising from $z$.

\section{Hochschild cohomology}
\subsection{}
The Hochschild cohomology $HH^*(\A,\A)$ of a (small) $A_\infty$-category $\A$
is a graded vector space, defined via the following cochain complex
$CC(\A,\A)$. A cochain of degree $r$ is a sequence $h = (h^d)_{d \geq 0}$ of
which each member $h^d$ is a family of graded linear maps
\[
\bigotimes_{i=1}^d hom_\A(X_i,X_{i+1}) \rightarrow hom_\A(X_1,X_{d+1})[r-d],
\]
one for each $(d+1)$-tuple $X_1,\dots,X_{d+1} \in \Ob\,\A$, satisfying
$h^d(a_d,\dots,a_1) = 0$ whenever some $a_i$ is the identity map $\id_{X_i}$.
In particular, $h^0$ consists of an element in $hom_\A^r(X,X)$ for each $X$.
The differential is
\[
\begin{split}
 & (\partial h)^d(a_d,a_{d-1},\dots,a_1) = \\
 & \qquad
 = \sum_{i+j \leq d+1}
 \mu_\A^{d+1-j}(a_d,\dots,a_{i+j},h^j(a_{i+j-1},\dots,a_i),
 a_{i-1},\dots,a_1)
 \\
 & \quad \qquad
 + \sum_{i+j \leq d+1}
 h^{d+1-j}(a_d,\dots,a_{i+j},\mu_\A^j(a_{i+j-1},\dots,a_i),
 a_{i-1},\dots,a_1).
\end{split}
\]
When $\A$ is directed, $CC(\A,\A)$ and hence $HH^*(\A,\A)$ are
finite-dimensional.

Hochschild cohomology has two related but different interpretations. One is in
terms of first order deformations of the maps $\mu^d_\A$, and $\partial h = 0$
appears there as the linearization of the structure equations of an
$A_\infty$-category, extended by an $h^0$ term. The other involves the
$A_\infty$-category $funct(\A,\A)$ of $A_\infty$-functors from $\A$ to itself,
and it says that
\begin{equation} \label{eq:functors}
HH^*(\A,\A) = H(hom_{funct(\A,\A)}(\mathrm{Id}_\A,\mathrm{Id}_\A))
\end{equation}
is the space of $A_\infty$-natural transformations from the identity functor to
itself. In the case where $\mu^d_\A = 0$ for $d \geq 3$, which means $\A$ is a
dg-category, there is yet another approach which is essentially equivalent to
\eqref{eq:functors} but more in line with classical homological algebra. Let
$\A\!-\!mod\!-\!\A$ be the dg-category of dg-functors from $\A^{opp} \times \A$
to chain complexes of vector spaces, and $D(\A\!-\!mod\!-\!\A)$ its derived
category, defined by inverting quasi-isomorphisms as in \cite[Chapter
10]{bernstein-lunts94}. Take the dg-functor $\Delta_{\A} \in
\Ob\,\A\!-\!mod\!-\!\A$ which assigns to $(X,Y) \in \Ob(\A^{opp} \times \A) =
\Ob\,\A \times \Ob\,\A$ the complex $hom_\A(X,Y)$ (to motivate the notation, we
should say that if $\A$ had only a single object $Z$, objects of
$\A\!-\!mod\!-\!\A$ would be just dg-bimodules over the dg-algebra
$hom_\A(Z,Z)$, and $\Delta_\A$ would be $hom_\A(Z,Z)$ considered as a bimodule
over itself). Then
\begin{equation} \label{eq:bimodule}
HH^*(\A,\A) = \Hom^*_{D(\A\!-\!mod\!-\!\A)}(\Delta_\A,\Delta_\A).
\end{equation}
The proof goes by observing that the Hochschild complex coincides with
$hom_{\A\!-\!mod\!-\!\A}(B_{\A}, \Delta_{\A})$ where $B_{\A}$, the bar
construction, is quasi-isomorphic to $\Delta_{\A}$ and $K$-projective in the
sense of \cite[Definition 10.12.2.1]{bernstein-lunts94}.

\subsection{}
An informal principle says that $HH^*(\A,\A)$ should be an invariant of
$D^b(\A)$. It is unknown whether this is rigorously true, but there is a weaker
result which is sufficient for our purpose.

\begin{prop} \label{th:hochschild}
Mutation of a directed $A_\infty$-category does not change its Hochschild
cohomology.
\end{prop}

To prove this one goes through the list of moves \cite[Definition
5.2]{seidel00}. If $F: \A \rightarrow \B$ is a quasi-isomorphism, there are
natural maps
\[
HH^*(\A,\A) \rightarrow H(hom_{funct(\A,\B)}(F,F)) \leftarrow HH^*(\B,\B),
\]
both of which can be shown to be isomorphisms by looking at the spectral
sequence associated to the length filtration. Next, changing the gradings of
the morphism groups $hom_\A(X^i,X^k)$ by some $\sigma_i-\sigma_k$ does not
affect the Hochschild complex. Now consider $\A \rightsquigarrow c\A$. One can
always find a directed dg-category ${\mathcal C}$ with a quasi-isomorphism
${\mathcal C} \rightarrow \A$, and then $c{\mathcal C}$ is quasi-isomorphic to
$c\A$. Hence it is sufficient to prove $HH^*(c{\mathcal C},c{\mathcal C}) \iso
HH^*({\mathcal C},{\mathcal C})$, which reduces the problem to the case of
dg-categories. In that case there is a derived Morita equivalence between
${\mathcal C}$ and $c{\mathcal C}$ , which means that there are objects $P \in
\Ob\,D(c{\mathcal C}\!-\!mod\!-\!{\mathcal C})$, $Q \in \Ob\,D({\mathcal
C}\!-\!mod\!-\!c{\mathcal C})$ with
\[
 P \otimes^{{\bf L}}_{{\mathcal C}} Q \iso \Delta_{c{\mathcal C}}, \quad
 Q \otimes^{{\bf L}}_{c{\mathcal C}} P \iso \Delta_{{\mathcal C}};
\]
for the definition of $\otimes^{{\bf L}}$ see again \cite[Chapter
10]{bernstein-lunts94}. The desired equality is a formal consequence of this
and \eqref{eq:bimodule}. The argument for $\A \rightsquigarrow r\A$ is
identical.

\begin{additional}
This section is essentially a review of known facts. $HH^*(\A,\A)$ is defined
in \cite{penkava-schwarz95}; the $A_\infty$-categories $funct(-,-)$ in
\cite{fukaya97}; and an account of derived Morita equivalences can be found in
\cite{koenig-zimmermann98}. These are just sample references, and by no means
the first ones historically.
\end{additional}

\section{Hochschild cohomology and global monodromy\label{sec:donaldson}}
\subsection{}
What follows is my attempt to formulate an idea of Donaldson. Let $(E,\pi)$ be
an exact Morse fibration over $D$ with a relative Maslov map $\delta_{E/D}$,
$M$ the fibre over the base point with its induced Maslov map $\delta_M$, and
$\mu \in \Sympe(M)$ the global monodromy, defined to be the symplectic parallel
transport around $\partial D$ in positive direction. This comes with a
canonical grading $\tilde{\mu}$, characterized by being zero near $\JJ_M \times
\partial M$. After making an admissible choice of paths $(c_1,\dots,c_m)$, with
corresponding distinguished basis $\Gamma = (L_1,\dots,L_m)$, one has that
$\tilde{\mu}$ is isotopic to $\ttau_{L_1}\dots\ttau_{L_m}$ within the group of
graded exact symplectic automorphisms. Choose gradings $\tL_i$ and let $\A =
\Lag(\widetilde{\Gamma})$ be the directed Fukaya category. By Proposition
\ref{th:hochschild}, $HH^*(\A,\A)$ is independent of the choice of $\Gamma$ and
of the gradings $\widetilde{\Gamma}$.

\begin{conjecture} \label{th:donaldson}
There is a long exact sequence
\begin{equation} \label{eq:donaldson}
\xymatrix{
 {HF^*(\tilde{\mu},+)} \ar[r] &
 {H^*(E;\Z/2)} \ar[d] \\
 & {HH^*(\A,\A)} \ar[ul]
}
\end{equation}
with the $\nwarrow$ map having degree one, and the others degree zero.
\end{conjecture}

It is elementary to check, using the classical Picard-Lefschetz formula, that
the Euler characteristics add up in the right way. There is also an informal
argument showing that, in a somewhat loose sense, $HF^*(\tilde{\mu},+)$ is
``constructed'' from nearly the same pieces as $HH^*(\A,\A)[-1]$, with
$H^*(E;\Z/2)$ compensating for the small difference. We will now reproduce this
argument, both because it seems to have been Donaldson's original motivation
and because it uses nicely the two exact sequences in Floer cohomology. As in
the discussion at the end of Section \ref{sec:tqft2}, it is convenient to use
the notation $HF(l,\pm)$ and $HF(c)$ for Floer cohomology groups.
\includefigure{skd1}{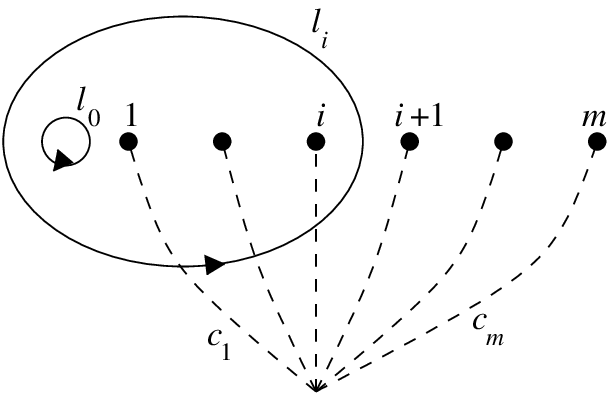}{hb}

\begin{myitemize}
\item[(i)]
Let $l_0,\dots,l_m$ be the loops in $D \setminus D^{\crit}$ shown in Figure
\ref{fig:skd1}. To each of them is associated a Floer group $HF(l_k,+)$. Our
aim is to successively ``decompose'' $HF(l_m,+)$. To begin, by Theorem
\ref{th:exact-sequence2} one has long exact sequences
\[
\xymatrix{
 {HF(l_k,+)} \ar[r] & {HF(l_{k-1},+)} \ar[r] &
 {HF(d_{1,k}),} \ar@/^1.5pc/[ll]^{ }
}
\]
where the $d_{1,k}$ are as in Figure \ref{fig:skd2}. We interpret this as
saying that $HF(l_m,+)$ is ``constructed'' by putting together the pieces
$HF(l_0,+)$, $HF(d_{1,1})$, \dots, $HF(d_{1,m})$.
\includefigure{skd2}{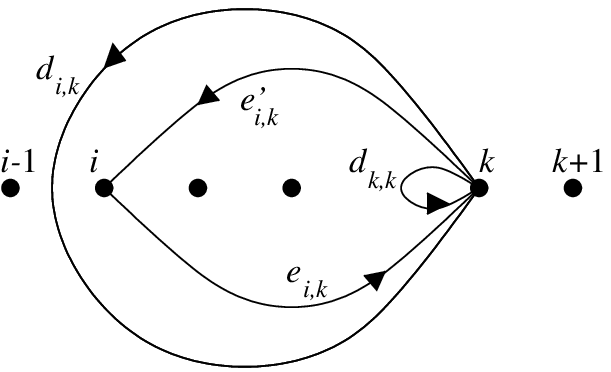}{hb}

\item[(ii)]
We next analyze $HF(d_{1,k})$ further. Take the paths $e_{i,k}$ and $e'_{i,k}$,
$i<k$, from Figure \ref{fig:skd2}. \cite[Theorem 3.3]{seidel00} yields long
exact sequences
\[
\xymatrix{
 {HF(d_{i,k})} \ar[r] & {HF(d_{i+1,k})} \ar[r] &
 {HF(e_{i,k}) \otimes HF(e'_{i,k}),} \ar@/^1.5pc/[ll]^{ }
}
\]
so that $HF(d_{1,k})$ is ``built'' from $HF(d_{k,k})$ and $HF(e_{i,k}) \otimes
HF(e'_{i,k})$, $1 \leq i < k$.

\item[(iii)]
A similar argument as in (ii) shows that $HF(e_{i,k}')$ can be ``decomposed''
into $HF(e_{i,k})^\vee$ together with $HF(e_{i,j})^\vee \otimes HF(e_{j,k}')$
for all $i < j < k$ (the dual vector space appears because the paths $e_{i,k}$,
$e_{i,j}$ occur with reversed orientation). One can apply the same argument
recursively to the second component $HF(e_{j,k}')$ of the tensor product, and
this result in a different ``decomposition'' of $HF(e_{i,k}')$, of which the
pieces are $(HF(e_{i_d,i_{d+1}}) \otimes \dots \otimes HF(e_{i_1,i_2}))^\vee$
ranging over all $d \geq 1$ and $i = i_1 < i_2 < \dots < i_{d+1} = k$.
\end{myitemize}
Putting all of this together yields a list of ``elementary pieces'' which make
up $HF(l_m,+)$. At this point we revert to the conventional notation, so that
$HF(l_m,+) = HF(\mu,+)$, $HF(l_0,+) = HF(\id_M,+) \iso H^*(M;\Z/2)$,
$HF(d_{k,k}) = HF(L_k,L_k) \iso H^*(L_k;\Z/2)$, and $HF(e_{i,k}) =
HF(L_i,L_k)$. We also restore the gradings, which were omitted up to now for
the sake of simplicity. Then the ``elementary pieces'' are:
\[
\left\{
\begin{split}
 & H^*(M;\Z/2); \\
 & (\Z/2)^m[-n\!-\!1]; \\
 & (\Z/2)^m[-1]; \text{ and} \\
 & HF^*(\tL_{i_1},\tL_{i_{d+1}})[-d\!-\!1] \otimes
 (HF^*(\tL_{i_d},\tL_{i_{d+1}}) \otimes \dots \otimes
 HF^*(\tL_{i_1},\tL_{i_2}))^\vee \\
 & \hspace{10em} \text{ for all $d \geq 1$ and $1\leq i_1 < \dots < i_{d+1} \leq m$.}
\end{split}
\right.
\]
The second and third piece are the result of splitting each $H^*(L_k;\Z/2)[-1]
= \Z/2[-1] \oplus \Z/2[-n\!-\!1]$. Now $E$ is obtained up to homotopy
equivalence by attaching $m$ cells of dimension $(n+1)$ to $M$. Hence the first
two pieces can be ``put together'' again, in the same informal sense as before,
to form $H^*(E;\Z/2)$. On the other hand the spectral sequence induced by the
length filtration, which converges to $HH^*(\A,\A)[-1]$, has $E_1$ term equal
to the direct sum of the remaining pieces (that is, all except the first two).
Of course this sort of bookkeeping falls far short of a proof of Conjecture
\ref{th:donaldson}, since it ignores all the maps in the long exact sequences,
and equally the differentials in the spectral sequence.

\subsection{}
For a better understanding of Conjecture \ref{th:donaldson} one would need to
know what the maps in it are. We outline here a possible description of one of
them, namely $H^*(E;\Z/2) \rightarrow HH^*(\A,\A)$. This is pure speculation,
supported mainly by the fact that the construction, apart from technicalities,
seems natural.

One may suppose that the $L_i$ are in generic position. Choose Morse functions
$f$ on $M$ and $f_i$ on $L_i$, such that $f|L_i = f_i$ and the restriction of
$f$ to a tubular neighbourhood of each $L_i$ is given by the sum of $f_i$ and a
positive definite quadratic form on the normal bundle to $L_i$. Moreover,
$f|\partial M = const.$ and the gradient should point outwards there. After
choosing suitable Riemannian metrics, one gets Morse cohomology complexes
$C_{Morse}(M)$ and $C_{Morse}(L_i)$ with $\Z/2$ coefficients, together with
natural restriction maps $r_i: C_{Morse}(M) \rightarrow C_{Morse}(L_i)$.
Enlarge $C_{Morse}(L_i)$ to a complex $\tilde{C}_{Morse}(L_i)$ by adding one
generator in degree $-1$ whose boundary is the only zero-dimensional cocycle
(the sum of all local minima of $f_i$). Composing $r_i$ and the inclusions
$C_{Morse}(L_i) \hookrightarrow \tilde{C}_{Morse}(L_i)$ gives a map $\tilde{r}:
C_{Morse}(M) \rightarrow \bigoplus_i \tilde{C}_{Morse}(L_i)$, whose mapping
cone satisfies
\begin{equation} \label{eq:morse-homology}
H^*(\Cone(\tilde{r})) \iso H^*(E;\Z/2)[1].
\end{equation}
The point of using this rather strange model for the cohomology of $E$ is that
there are natural maps
\begin{equation} \label{eq:first-order}
\tilde{C}_{Morse}(L_i) \rightarrow CC(\A,\A)[1], \quad C_{Morse}(M) \rightarrow
CC(\A,\A).
\end{equation}
The first of these is a chain homomorphism; the second isn't, but the direct
sum of the two is a chain map $\Cone(\tilde{r}) \rightarrow CC(\A,\A)[1]$,
which in combination with \eqref{eq:morse-homology} gives the desired map
$H^*(E;\Z/2) \rightarrow HH^*(\A,\A)$.

Both maps in \eqref{eq:first-order} are constructed by the method of ``cutting
down moduli spaces'' which is familiar from the definition of the cap product
on Floer cohomology. We will describe the first map in detail, and the second
one only briefly. Fix $i \in \{1,\dots,m\}$ and let $Z \subset L_i$ be the
unstable manifold, under the gradient flow, of some critical point of $f_i$ of
Morse index $p$. Suppose for the moment that $p>0$, so that $Z$ has positive
codimension. Recall that the boundary operator $\mu_{\A}^1$ on $CF^*(\tL_{i_1},
\tL_{i_2})$, $i_1<i_2$, is defined using pseudo-holomorphic maps $u: \R \times
[0;1] \rightarrow M$ with boundary conditions $u(\R \times \{2-\nu\}) \subset
L_{i_\nu}$. If $i_\nu = i$ for some $\nu \in \{1,2\}$, one can consider the
subset of those maps $u$ which satisfy $u(0,2-\nu) \in Z$. If one assumes that
the almost complex structures have been chosen generically (similar assumptions
will be made tacitly later on) then counting isolated points in that subset
yields a map
\begin{equation} \label{eq:h1}
h^1_Z: CF^*(\tL_{i_1},\tL_{i_2}) \rightarrow CF^{*+p}(\tL_{i_1},\tL_{i_2}).
\end{equation}
Extend this to the other Floer cochain groups in $\A$, that is to say to
$CF^*(\tL_{i_1},\tL_{i_2})$ with $i_1,i_2 \neq i$, by setting $h_Z$ to be zero
there. Now consider the composition maps of order $d \geq 2$, $\mu_\A^d:
CF^*(\tL_{i_d},\tL_{i_{d+1}}) \otimes \dots \otimes CF^*(\tL_{i_1}, \tL_{i_2})
\rightarrow CF^{*+2-d}(\tL_{i_1},\tL_{i_{d+1}})$ for $i_1 < \dots < i_{d+1}$.
Their definition uses moduli spaces of pairs $(r,u)$, where $r \in \RR^{d+1}$
is a point in the moduli space of discs with $d+1$ marked boundary point,
$\SS^{d+1,*}_r$ is the corresponding disc with the marked points removed, and
$u: \SS^{d+1,*}_r \rightarrow M$ is a map satisfying a pseudo-holomorphicity
equation. $\partial \SS^{d+1,*}_r$ consists of intervals $I^{d+1}_{r,\nu}$, $1
\leq \nu \leq d+1$, and the boundary conditions are
\[
u(I^{d+1}_{r,\nu}) \subset L_{i_\nu}.
\]
Assuming that $i_{\mu} = i$ for some $\mu$, one considers the moduli space of
triples $(r,u,z)$, where $z$ is a point in $I^{d+1}_{r,\mu}$ such that $u(z)
\in Z$. Counting isolated points in this moduli space defines a map
\[
h^d_Z: CF^*(\tL_{i_d},\tL_{i_{d+1}}) \otimes \dots \otimes CF^*(\tL_{i_1},
\tL_{i_2}) \rightarrow CF^{*+p+1-d}(\tL_{i_1},\tL_{i_{d+1}}).
\]
Again, this is extended trivially to the cases when $i_1,\dots,i_{d+1} \neq i$.
The sequence $(h_Z^0 = 0,h_Z^1,h_Z^2,\dots)$ is a Hochschild cochain $h_Z \in
CC^{p+1}(\A,\A)$.

The case $p = 0$, which we had excluded above, is considerably easier: one
takes \eqref{eq:h1} to be the projection to the subspace spanned by the points
of $L_{i_1} \cap L_{i_2}$ which lie on $Z$, and all other $h^d_Z$ to be zero.
Finally there is a canonical element $\tilde{h} \in CC^0(\A,\A)$, defined by
setting $\tilde{h}^0$ to be the unit element in $hom_\A^0(\tL_i,\tL_i) \iso
\Z/2$, with all other components equal to zero. The assignment $Z \mapsto h_Z$,
extended to $\tilde{C}_{Morse}(L_i)$ by mapping the generator of degree $-1$ to
$\tilde{h}$, gives the first map in \eqref{eq:first-order}.

The definition of the second map uses moduli spaces of triples $(r,u,z)$ with a
point $z$ which can lie anywhere in $\SS^{d+1,*}_r$, such that $u(z)$ is in
some stable manifold of $f$. The failure of this to be a chain map reflects the
fact that $z$ can move to the boundary; in that case however $u(z)$ necessarily
lies on some $L_i$, which provides the connection with the first map.

\begin{remarks}
(i) This discussion may seem rather abstract, but it has concrete computational
implications, most strikingly for $n = 1$ when $M$ is a surface. In that case
the $A_\infty$-category $\A$ is constructed from $\Gamma$ by a purely
combinatorial count of immersed polygons, and the same is true for the maps
\eqref{eq:first-order}. Suppose that Conjecture \ref{th:donaldson} is true and
that our guess for the map $H^*(E;\Z/2) \rightarrow HH^*(\A,\A)$ which occurs
in it is correct. Then one gets a combinatorial algorithm for computing the
Floer cohomology $HF(\tau_{L_1}\dots\tau_{L_m},+)$ of an arbitrary product of
Dehn twists in $M$ (along exact curves $L_i$ which admit gradings, to be
precise).

(ii) From a deformation theory point of view, what we have described is a first
order infinitesimal deformation of $\A$ parame\-trized by the graded vector
space $H^*(E;\Z/2)$. This was done by inserting one marked point into the
Riemann surfaces which define $\mu^d_\A$. It seems natural to expect that using
more marked points will allow one to extend this to higher order. This is
interesting because in principle there are obstructions to such an extension,
which are expressed by a natural graded Lie bracket on $HH^*(\A,\A)[1]$. If the
extension to higher order is indeed possible, it would mean that the Lie
bracket vanishes on the image of $H^*(E;\Z/2) \rightarrow HH^*(\A,\A)$.
\end{remarks}

\subsection{}
To round off the discussion of Conjecture \ref{th:donaldson} we look at its
implications in some specific cases. Suppose first that $(E,\pi)$ and
$\delta_{E/D}$ come from a Lefschetz pencil satisfying
\eqref{eq:relative-maslov} for some $a \in \Z$. Then $\mu$ can be isotoped to
the identity but only within a group of symplectic automorphisms which is
larger than $\Sympe(M)$, the isotopy being nontrivial along $\partial M$. Floer
cohomology is not invariant under such isotopies. However, a more careful
analysis shows that $HF(\mu,+) \iso HF(\id_M,-) \iso H^*(M,\partial M;\Z/2)$.
More precisely, taking the canonical grading into account, this formula reads
\[
HF^*(\tilde{\mu},+) \iso H^{*+a-4}(M,\partial M;\Z/2).
\]
For $X = \CP{2}$, $\xi = {\mathcal O}(2)$, $\A$ has been determined in Section
\ref{sec:examples}a, and one finds by explicit computation that $HH^r(\A,\A)
\iso H^r(M,\partial M;\Z/2) \oplus H^r(E;\Z/2)$ for all $r$. This is compatible
with an exact sequence \eqref{eq:donaldson} in which the horizontal arrow would
be the zero map. There are reasons to believe that this arrow will in fact
vanish for all Lefschetz pencils except the trivial one (the degree one pencil
on $X = \CP{n+1}$).
\includefigure{mirror-monodromy}{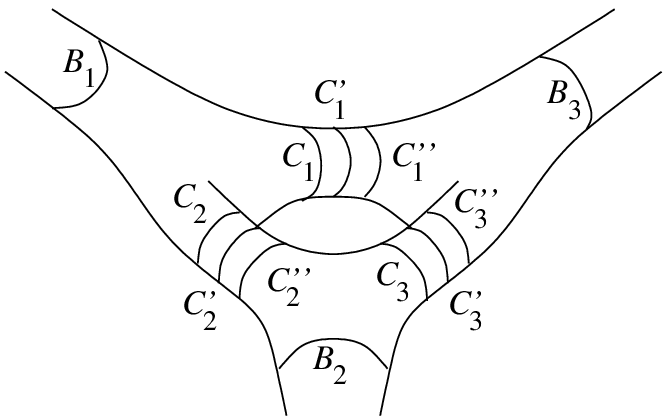}{hb}

Next consider the mirror dual of the projective plane, as in Section
\ref{sec:examples}b. From the structure of the singular fibre $Y_\infty$ one
sees that the global monodromy can be written as a product of Dehn twists and
their inverses along disjoint curves, namely
\[
\mu = \prod_{i=1}^3
(\tau_{B_i}\tau_{C_i}^{-1}\tau_{C_i'}^{-1}\tau_{C_i''}^{-1})
\]
where the curves are as shown in Figure \ref{fig:mirror-monodromy}. Using
\cite{seidel96b} one computes that $HF^*(\tilde{\mu},+) = H^{*-2}(M \setminus
\bigcup_{i=1}^3 C_i \cup C_i' \cup C_i'',\partial M;\Z/2)$. Comparing
dimensions shows that $HH^r(\A,\A) \iso HF^{r+1}(\tilde{\mu},+) \oplus
H^r(E;\Z/2)$, as in the previous example.

As a third and final case, suppose that $(E,\pi)$ and $\delta_{E/D}$ come from
a Morsification of an isolated hypersurface singularity. We conjecture that for
all nontrivial singularities
\begin{equation} \label{eq:acyclicity}
HF(\mu) = 0;
\end{equation}
this can be checked easily in many cases, e.g.\ for weighted homogeneous
singularities. From \eqref{eq:acyclicity} and Conjecture \ref{th:donaldson} it
would follow that $HH^*(\A,\A) \iso H^*(E;\Z/2) = \Z/2$. Because of the role of
Hochschild cohomology in deformation theory, this means that $\A$ is rigid: it
admits no nontrivial first order deformations of any degree. As of now, there
is no deeper explanation for this phenomenon.

\begin{additional}
For the cap product on Floer cohomology see \cite{le-ono97}, \cite{schwarz96}.
Equation \eqref{eq:acyclicity}, if true, improves on old topological results of
A'Campo \cite{acampo73b} and L{\^e} \cite{le74}.
\end{additional}

\section{Morse categories\label{sec:morse}}
\subsection{}
Let $N^{n+1}$ be a compact manifold with codimension two corners. More
precisely, we assume that $\partial N$ consists of three codimension one faces
$\partial_- N$, $\partial_0 N$, $\partial_+ N$ (any of which may be empty) with
$\partial_-N \cap \partial_+N = \emptyset$, and where the corners are
$(\partial_- N \cup \partial_+N) \cap \partial_0 N$. Let $p: N \rightarrow
[-1;1]$ be a Morse function satisfying $\partial_-N = p^{-1}(-1)$, $\partial_+N
= p^{-1}(1)$, and such that $p|\partial_0N$, as a function on $\partial_0N$,
has no critical points. In addition we require that if $x \in N$ is a critical
point of $p$ with Morse index $i(x)$, then
\begin{equation} \label{eq:self-indexing}
 \begin{cases}
 -1 < p(x) < 0 & \text{if $i(x) = 0$,} \\
 p(x) = 0 & \text{if $0<i(x)<n+1$,} \\
 0 < p(x) < 1 & \text{if $i(x) = n+1$.} \end{cases}
\end{equation}
For $n = 1$ this is more or less the classical notion of self-indexing Morse
function, but in higher dimensions it is far more restrictive. To any function
$p$ with these properties one can associate a {\bf Morse category}
$Morse(N,p)$. To do that, choose a Riemannian metric such that $\nabla p$ is
parallel to $\partial_0N$. For two critical points $x,y$ of $p$, let ${\mathcal
G}(x,y)$ be the space of unparametrized gradient trajectories going from $x$ to
$y$. This is a smooth manifold; it is compact unless $i(x) = 0$ and $i(y) =
n+1$, in which case there is a natural compactification to a manifold with
boundary $\overline{\mathcal G}(x,y)$ such that
\begin{equation} \label{eq:compactification}
\partial\overline{\mathcal G}(x,y) = \bigsqcup_{0 < i(w) < n+1} {\mathcal
G}(w,y) \times {\mathcal G}(x,w).
\end{equation}
These are all standard facts from Morse theory. What is important is that there
is no need to choose the metric generically, since \eqref{eq:self-indexing}
implies that $\nabla p$ is always Morse-Smale. A cobordism argument based on
this fact shows that the topological type of the spaces ${\mathcal G}(x,y)$ and
of their compactifications is independent of the metric.

The objects of $\A = Morse(N,p)$ are the critical points of $p$, and the
morphisms are $\hom_\A(x,y) = H_{-*}({\mathcal G}(x,y);\Z/2)$; note the sign in
front of the grading. The composition maps
\begin{equation} \label{eq:multiplication}
H_{-*}({\mathcal G}(w,y)) \otimes H_{-*}({\mathcal G}(x,w)) \rightarrow
H_{-*}({\mathcal G}(x,y))
\end{equation}
are as follows. If $x = w$ then $[point] \in H_0(point;\Z/2) = hom^0_\A(x,w)$
acts as identity morphism, and similarly for $w = y$. If $i(x) = 0$, $0 < i(w)
< n+1$, and $i(y) = n+1$, the composition comes from the inclusion of the
boundary stratum \eqref{eq:compactification}, together with the fact that
${\mathcal G}(x,y) \hookrightarrow \overline{\mathcal G}(x,y)$ is a homotopy
equivalence. In all remaining cases \eqref{eq:multiplication} is automatically
zero, because the tensor product on the left hand side vanishes. Thus $\A$ is a
$\Z/2$-linear and $\Z$-graded category. Ordering the critical points according
to $p(x)$ shows that it can also be viewed as a directed $A_\infty$-category
with $\mu^d_\A = 0$ for $d \neq 2$.

Let $B(N,p) \subset N \setminus \partial N$ be the subspace of points whose
gradient flow line exists for all time, which means that it converges to a
critical point in both positive and negative direction. Define the fundamental
object $\fund{B(N,p)} \in \Ob\,D^b(\A)$ to be the twisted complex
\[
(C,\delta_C) = \Big(\bigoplus_{dp(x) = 0} x[-i(x)], (\delta_{xy}) \Big),
\]
where $\delta_{xy} \in hom^1_\A(x[-i(x)],y[-i(y)])$ is the fundamental homology
class $[{\mathcal G}(x,y)]$ if that space is compact and $x \neq y$, and zero
otherwise. The generalized Maurer-Cartan equation reduces to
$\mu^2_\A(\delta_C,\delta_C) = 0$, which is a consequence of
\eqref{eq:compactification}. The notation $\fund{B(N,p)}$ is motivated by the
fact that
\begin{equation} \label{eq:sheafcohomology}
\Hom^*_{D^b(\A)}(\fund{B(N,p)},\fund{B(N,p)}) \iso H^*(B(N,p);\Z/2).
\end{equation}
To see this one identifies $hom_{\Tw\,\A}(C,C)$ with the $E^1$ term of the
spectral sequence, converging to $H^*(B(N,p);\Z/2)$, which arises from a
decomposition of $B(N,p)$ determined by the gradient flow. There is only one
nonzero differential in this spectral sequence, and that coincides with the
differential on $hom_{\Tw,\A}(C,C)$, which completes the proof. If $N$ is
closed, so that $B(N,p) = N$, the fundamental object $\fund{N}$ has another
property. Namely, for any $X \in \Ob\,D^b(\A)$ the composition gives rise to a
nondegenerate pairing
\begin{equation} \label{eq:verdier}
 \Hom^*_{D^b(\A)}(X,\fund{N}) \otimes \Hom^{n+1-*}_{D^b(\A)}(\fund{N},X)
 \rightarrow
 \Hom^{n+1}_{D^b(\A)}(\fund{N},\fund{N}) \iso \Z/2.
\end{equation}

\begin{remarks}
(i) Consider $B(N,p)$ as a space stratified by the stable manifolds of $\nabla
p$. It was pointed out to me by Khovanov that the bounded derived category of
sheaves of $\Z/2$-vector spaces, constructible with respect to this
stratification, is equivalent to $D^b(\A)$; this appears to be a reformulation
of a familiar result \cite{kapranov90}. The objects $x[-i(x)]$ correspond to
the constant sheaves along the strata, and $\fund{B(N,p)}$ to the constant
sheaf on all of $B(N,p)$. This explains the properties listed above:
\eqref{eq:sheafcohomology} computes the cohomology of the constant sheaf, and
\eqref{eq:verdier} is Verdier duality.

(ii) The Morse categories for $p$ and $-p$ are not equivalent. However, they
are Koszul dual, so that their derived categories are equivalent.
\end{remarks}

\subsection{}
Assume now that $n = 1$, so that $p^{-1}(0)$ is a nodal curve on the surface
$N$. Put a sign $\ominus,\oplus$ into each connected component of $N \setminus
p^{-1}(0)$ that contains a local minimum respectively maximum of $p$. This
diagram retains enough information about $p$ to reconstruct $\Morse(N,p)$.
\includefigure{triple-point}{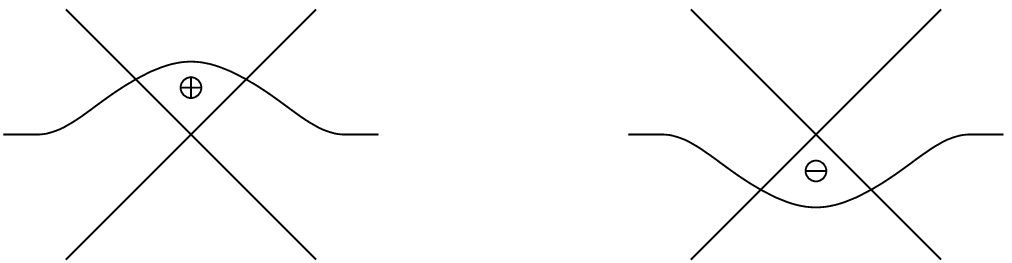}{ht}%

Let's say that two functions are related by a triple point move if the
corresponding diagrams differ as shown in Figure \ref{fig:triple-point}. By a
sequence of mutations that mimics the base changes in \cite[Volume\ 2, I \S 4]
{arnold-gusein-zade-varchenko}, one can prove

\begin{prop} \label{th:triple-point}
Up to equivalence, a triple point move does not change the derived Morse
category $D^bMorse(N,p)$.
\end{prop}

It would be interesting to find a topological interpretation of $HH^*(\A,\A)$
for $\A = Morse(N,p)$. I have not succeeded in doing that, except for one
special situation. Call $(N,p)$ cellular if all spaces ${\mathcal G}(x,y)$ are
contractible (we are still assuming that $n = 1$; anyway cellularity is
impossible in higher dimensions, except for trivial cases).

\begin{prop} \label{th:cellular}
If $(N,p)$ is cellular, $HH^*(\A,\A) \iso H^*(B(N,p);\Z/2)$.
\end{prop}

The first step in the proof is to note that the only possible nonzero groups
are $HH^r(\A,\A)$ for $r = 0,1,2$. $HH^0(\A,\A) \iso H^0(B(N,p);\Z/2)$ is
always true. From the cellularity one derives that $HH^1(\A,\A) \iso
H^1(B(N,p);\Z/2)$, and also that the Euler characteristic of $HH^*(\A,\A)$
equals that of $B(N,p)$, which completes the argument.

\begin{additional}
Objects similar to our Morse categories were introduced several years ago by
Cohen, Jones and Segal \cite{cohen-jones-segal95} under the name flow
categories. Proposition \ref{th:cellular} was inspired by A'Campo's paper
\cite{acampo00}.
\end{additional}

\section{Real structures\label{sec:real}}
\subsection{}
For the duration of this section, we extend the definition of exact Morse
fibration by allowing each fibre to contain several critical points. This does
not affect the theory seriously. The notion of admissible choice of paths
remains the same; each path may give rise to several disjoint vanishing cycles,
which one can place in arbitrary order in the distinguished basis. The
important thing is that such bases are still unique up to Hurwitz moves.

Let $(E,\pi) = (E,\pi,\Omega,\Theta,J_0,j_0)$ be an exact Morse fibration over
$D$, with $j_0$ the standard complex structure on $D$, and $\delta_{E/D}$ a
relative Maslov map. A real structure on it is an involution $\iota: E
\rightarrow E$ such that
\[
 \pi(\iota(x)) = \overline{\pi(x)},\;\; \iota^*\Omega = -\Omega,\;\;
 \iota^*\Theta = -\Theta,\;\; \iota^*J_0 = -J_0.
\]
There is always a $(j,J) \in \JJ_{E/D}$ with $j$ the standard complex structure
on $D$ and such that $\iota^*J = -J$. Then the line bundle
\begin{equation} \label{eq:delta}
\Delta_{E/D} \,|\, \{(j,J)\} \times E = \Lambda^{n+1}(TE,J)^{\otimes 2} \otimes
\pi^*(TD,j)^{\otimes -2}
\end{equation}
has a natural real involution $\iota_{(j,J)}$ which covers $\iota$. Call
$\iota$ compatible with the relative Maslov map if $\delta_{E/D}\,|\,\{(j,J)\}
\times E$ is homotopic to a trivialization of \eqref{eq:delta} which takes
$\iota_{(j,J)}$ to complex conjugation in $\C$. This condition is independent
of the choice of $(j,J)$.

Given a real structure, one can consider $N = E^\iota$ with the Morse function
$p = \pi|N: N \rightarrow [-1;1]$. $N$ is a manifold with corners of the kind
considered in the previous section, with $\partial_-N \cup \partial_+N = (N
\cap \partial_vE) = p^{-1}(\{\pm 1\})$ and $\partial_0N = N \cap \partial_hE$;
and $p|\partial_0N$ has no critical points because the same holds for
$\pi|\partial_hE$. We say that $\iota$ a good real structure if the critical
points of $p$ satisfy \eqref{eq:self-indexing}, in which case there is a
well-defined Morse category $Morse(N,p)$. We say that $\iota$ is complete if
all critical points of $\pi$ lie on $N$.

\begin{conjecture} \label{th:real}
Let $(E,\pi)$ be an exact Morse fibration over $D$ with a good and complete
real structure $\iota$, which is also compatible with a relative Maslov map
$\delta_{E/D}$. Take an admissible choice of paths $(c_1,\dots,c_m)$, each of
which lies in the lower half-disc $D \cap \{\im(z) \leq 0\}$, and let $\Gamma$
be the corresponding distinguished basis. Then $\Lag(\widetilde{\Gamma})$, for
a certain choice of grading $\widetilde{\Gamma}$, is quasi-isomorphic to
$Morse(N,p)$.
\end{conjecture}

For $n = 1$ there is an elementary way of drawing the vanishing cycles which
constitute $\Gamma$, starting from $(N,p)$; see e.g.\ \cite{acampo98}. Using it
I have proved Conjecture \ref{th:real} in that dimension. There is also a
general argument in support of it which, while falling short of a proof, is
quite suggestive. Take an interval $[a;b] \subset [-1;1] = D \cap \R$ such that
$[a;b] \cap D^\crit = \{a;b\}$. The symplectic parallel transport along $(a;b)$
is $\iota$-equivariant, and hence preserves the real part $N$; in fact its
restriction to $N$ is equal, up to reparametrization of the orbits, to the flow
of $\nabla p$ with respect to some metric. Using property
\eqref{eq:self-indexing} one can show that intersection points of the vanishing
cycles arising from the endpoints $\{a,b\}$ are in one-to-one correspondence
with gradient flow lines joining the corresponding critical points of $p$.

\begin{remark}
The combination of Conjecture \ref{th:real}, Proposition \ref{th:cellular} and
Conjecture \ref{th:donaldson} leads to intriguing relations between the global
monodromy of an exact Morse fibration on one hand, and its possible good and
complete real structures on the other hand.
\end{remark}

Real structures are not the most general context in which there is a relation
between directed Fukaya categories and Morse categories. Take an exact Morse
fibration $(E,\pi)$ over $D$ with relative Maslov map $\delta_{E/D}$. A
pseudo-real submanifold is a pair of embeddings
\begin{equation} \label{eq:pseudo-real}
\xymatrix{
 {N} \ar[d]_{p} \ar[r]^-{C} & E \ar[d]^{\pi} \\
 {[-1;1]} \ar[r]^{c} & D.
}
\end{equation}
Here $c$ is a curve in $D$ with $c^{-1}(\partial D) = \{\pm 1\}$, meeting
$\partial D$ transversally. $N$ is an $(n+1)$-dimensional manifold with corners
as before. $C$ is an embedding transverse to $\partial_hE$ and $\partial_vE$,
such that $C^{-1}(\partial_vE) = \partial_-N \cup \partial_+N$ and
$C^{-1}(\partial_hE) = \partial_0N$; we require that $C^*\Omega$ is zero and
$C^*\Theta$ is an exact one-form. $p$ is a Morse function such that
$\partial_{\pm}N = p^{-1}(\pm 1)$. Furthermore, using the fact that $\im(C)$ is
totally real with respect to $J$ for any $(j,J) \in \JJ_{E/D}$, one finds that
there is a canonical homotopy class of sections of $S(\Delta_{E/D})$ over
$\{(j,J)\} \times \im(C)$. If a section in that homotopy class can be chosen
such that its composition with $\delta_{E/D}$ is a constant map $N \rightarrow
S^1$, we say that our pseudo-real submanifold is compatible with the relative
Maslov map. The notion of good pseudo-real submanifold is defined as before.

\begin{conjecture}
Suppose that $(N,p)$ is a good pseudo-real submanifold in $(E,\pi)$, compatible
with $\delta_{E/D}$. Then $D^b\Lag(\widetilde{\Gamma})$ contains
$D^bMorse(N,p)$ as a full triangulated subcategory. In particular there is an
object $\fund{B(N,p)} \in \Ob\,D^b\Lag(\widetilde{\Gamma})$ whose endomorphism
ring is $H^*(B(N,p);\Z/2)$.
\end{conjecture}

Again, it is not difficult to verify this for $n = 1$.

\subsection{}
As examples one can consider Lefschetz pencils as in Section
\ref{sec:examples}a which are real, meaning that $X$ and $\xi$ have real
involutions such that the composition of the two preserves $\sigma_0,\sigma_1$.
Then the resulting exact Morse fibration has a natural real structure. To get
the right kind of relative Maslov map, one must ask that the isomorphism
\eqref{eq:relative-maslov} be compatible with the involutions on both sides. An
example is $X = \CP{2}$, $\xi = {\mathcal O}(2)$, $\sigma_0(x) = x_0^2-x_2^2$,
$\sigma_1(x) = x_0^2 + x_1^2 + x_2^2$, with complex conjugation. This gives
rise to a complete and good real structure such that $(N,p)$ is $\RP{2}$ with
its standard Morse function (having three critical points). One verifies easily
that $Morse(N,p)$ is precisely given by the quiver with relations
\eqref{eq:cp2-quiver}. Another possible application is to isolated hypersurface
singularities which have real Morsifications. For instance, the Morsifications
used in Section \ref{sec:examples}c for ADE plane curve singularities are real,
and this provides another way of computing the associated categories.

\begin{additional}
The basic idea of relating vanishing cycles and real structures comes from work
of A'Campo \cite{acampo75b} and Gusein-Zade \cite{gusein-zade74}.
\end{additional}

\section{Matching pairs and matching paths\label{sec:induction}}
\subsection{}
Let ${\mathcal C}$ be a triangulated category, linear over $\Z/2$ and such that
the spaces $\Hom^*_{\mathcal C}(X,Y)$ are finite-dimensional. Suppose that
$(Y^1,\dots,Y^m)$ is a full exceptional collection in ${\mathcal C}$ such that
for some $1 \leq i < m$, the following properties hold. Firstly
\begin{equation} \label{eq:matching}
\Hom^*_{\mathcal C}(Y^i,Y^{i+1}) = H^*(S^n;\Z/2)
\end{equation}
for some $n>0$. Secondly, if $a$ denotes the unique nontrivial morphism $Y^i
\rightarrow Y^{i+1}$ of degree zero, the compositions
\begin{align*}
 & a \circ -: \Hom^*_{\mathcal C}(Y^k,Y^i) \longrightarrow
 \Hom^*_{\mathcal C}(Y^k,Y^{i+1}) && \text{for $k<i$, and} \\
 & - \circ\, a: \Hom^*_{\mathcal C}(Y^{i+1},Y^l) \longrightarrow
 \Hom^*_{\mathcal C}(Y^i,Y^l) && \text{for $l>i+1$}
\end{align*}
are isomorphisms. We then call $(Y^i,Y^{i+1})$ a {\bf matching pair} of
dimension $n$. Let $C \in \Ob\, {\mathcal C}$ be the cone of $a$, which is
well-defined up to isomorphism.

\begin{lemma} \label{th:matching-pair}
$C$ is a spherical object of dimension $(n+1)$ in the sense of Definition
\ref{def:spherical}.
\end{lemma}

For future use we need to mention one more fact about $C$. Assume that
${\mathcal C} = D^b(\A)$ is the derived category of some directed
$A_\infty$-category $\A$, with $\Ob\,\A = \{X^1,\dots,X^m\}$, and moreover that
$(Y^1,\dots,Y^m)$ is obtained from $(X^1,\dots,X^m)$ by mutation, that is to
say by applying the transformations \cite[Equations (6) and (7)]{seidel00} and
their inverses, plus possibly shifting each object by some amount. In this
situation the following holds:

\begin{lemma} \label{th:invariant}
If $F: \A \rightarrow \A$ is a quasi-isomorphism which preserves the numbering
of the objects, then $(D^bF)(C) \iso C$.
\end{lemma}

The proof is a general nonsense argument. Consider first a general triangulated
category ${\mathcal C}$ as before with an exact self-equivalence $G$, and
define $Inv(G) \subset \Ob\,{\mathcal C}$ to be the class of objects $X$ such
that $G(X) \iso X$. This has the following properties: (i) if $X$ lies in
$Inv(G)$ then so does $X[\sigma]$ for any $\sigma \in \Z$, as well as any
object isomorphic to $X$; (ii) if $X,Y \in Inv(G)$ and there is a unique
nonzero map $X \rightarrow Y$ of degree zero, then the cone over it lies in
$Inv(G)$; (iii) if $X,Y \in Inv(G)$ then $T_X(Y)$, $T_X'(Y) \in Inv(G)$. Here
$T_X(Y)$ is as in \cite[Section 5]{seidel00}, and $T_X'(Y)$ is the object which
fits into an exact triangle
\[
T_X'(Y) \rightarrow Y \xrightarrow{ev^\vee} \Hom^*_{\mathcal C}(Y,X)^\vee
\otimes X \rightarrow T_X'(Y)[1],
\]
where ${ev}^\vee$ is the transpose of the evaluation $ev: \Hom^*_{\mathcal C}
(Y,X) \otimes Y \rightarrow X$. Now specialize to ${\mathcal C} = D^b(\A)$, $G
= D^bF$. Then $X^i \in Inv(G)$, and by (iii) the same holds for the objects
$Y^i$ of any mutated exceptional collection. Applying (ii) shows that $C \in
Inv(G)$.

\subsection{}
Let $(E,\pi) = (E,\pi,\Omega,\Theta,J_0,j_0)$ be an exact Morse fibration over
a compact base $S$, with $\dim\,E = 2n+2$. As explained in \cite[Section
3]{seidel00}, to any smooth path $c: [0;1] \rightarrow S$ with $c^{-1}(S^\crit)
= \{0;1\}$ and $c'(0),c'(1) \neq 0$ one can associate two vanishing cycles
$V_{c,0},V_{c,1}$, which are exact framed Lagrangian spheres in, say,
$E_{c(1/2)}$. Now suppose that $c$ is an embedding, does not intersect
$\partial S$, and that its vanishing cycles are isotopic, which means that
there is a smooth family $I = (I_t)_{0 \leq t \leq 1}$ inside the space of all
exact Lagrangian submanifolds of $E_{c(1/2)}$ joining $I_0 = V_{c,0}$ and $I_1
= V_{c,1}$. Pairs $(c,[I])$ of this kind, where $[I]$ denotes the homotopy
class of $I$ rel endpoints, are called {\bf matching paths}. The interest of
this notion lies in a construction which I learnt from Donaldson. Roughly
speaking (see below for a precise formulation) this produces from a matching
path a submanifold of the total space $E$, which is Lagrangian in a suitable
sense. The idea is easy to see in the case where $V_{c,0} = V_{c,1}$ and $I$ is
the constant isotopy. By the definition of vanishing cycle, there are canonical
$(n+1)$-dimensional balls $D_{c,\nu} \subset E$, $\nu = 0,1$, such that
$\pi(D_{c,0}) = c([0;1/2])$, $\pi(D_{c,1}) = c([1/2;1])$, $\partial D_{c,\nu} =
V_{c,\nu}$ and $\Omega|D_{c,\nu} = 0$. Under the assumption that $V_{c,0} =
V_{c,1}$, the two balls fit together smoothly to give a closed submanifold of
$E$. Fix a positive two-form $\beta$ on $S$. Then the base $\Omega^{(r)} =
\Omega + r\cdot\pi^*\beta$ is symplectic for large positive $r$, and $D_{c,0}
\cup D_{c,1}$ is Lagrangian with respect to it.

In general, rather than choosing some $r$, it seems more natural to work with
the entire family of forms $\Omega^{(r)}$. This requires a suitably adapted
notion of Lagrangian submanifold. Let $\Lambda \subset \Rgeq \times (E
\setminus \partial E)$ be a $(n+2)$-dimensional submanifold such that the
projection $\Lambda \rightarrow \Rgeq$ is a proper submersion; this is the same
as a smooth family of closed $(n+1)$-submanifolds $\Lambda^{(r)} \subset E$,
$\Lambda^{(r)} = \Lambda \cap (\{r\} \times E)$. For simplicity assume that
$H^1(\Lambda;\R) = 0$. Such a $\Lambda$ is called eventually Lagrangian if
$\Omega^{(r)}|\Lambda^{(r)} = 0$ for $r \gg 0$. We can now state precisely the
nature of the construction mentioned above: for any matching pair $(c,[I])$ it
gives an eventually Lagrangian submanifold $\Lambda_{(c,[I])}$, unique up to
isotopy.

The first step is to trivialize $E$ over $c$, away from its endpoints, using
symplectic parallel transport. The trivialization is the unique embedding
$\Psi: (0;1) \times E_{c(1/2)} \rightarrow E$ which satisfies $\Psi(\{t\}
\times E_{c(1/2)}) = E_{c(t)}$, $\Psi\,|\,\{1/2\} \times E_{c(1/2)} = \id$, and
$\Psi^*\Omega = \Omega|E_{c(1/2)}$. Because $D_{c,0}$ and $D_{c,1}$ are
likewise defined by parallel transport, they satisfy $\Psi^{-1}(D_{c,0}) =
(0;1/2] \times V_{c,0}$, $\Psi^{-1}(D_{c,1}) = [1/2;1) \times V_{c,1}$. One can
assume that the isotopy $(I_t)$ is constant near the endpoints. Then
\[
 A = \{x_0\} \cup \{x_1\} \cup
     \bigcup_{t \in [0;1]} \Psi(\{t\} \times I_t),
\]
where $x_\nu \in E_{c(\nu)}$ are the critical points of $\pi$, is a closed
submanifold of $E \setminus \partial E$ and contained in $\pi^{-1}c([0;1])$.

By dualizing the normal vector field $\partial I_t/\partial t$ one gets a
function $h_t$ (unique up to a constant) on each $I_t$. Extend the $h_t$ to
functions $H_t$ on $E_{1/2}$ supported in a neighbourhood of $I_t$; these may
be chosen such that $H_t = 0$ for $t$ close to $0$ or $1$. One can then find an
$\alpha \in \Omega^1(E)$ vanishing near $\partial E$, which satisfies
$\alpha\,|\,\ker(D\pi) = 0$ and $\Psi^*\alpha = H_t\,dt$. With respect to the
modified symplectic forms $\Omega^{(r)} + d\alpha$, $r \gg 0$, $A$ is a
Lagrangian submanifold. Moser's Lemma shows that $\Omega^{(r)} + d\alpha$ and
$\Omega^{(r)}$ are diffeomorphic, and the diffeomorphism can be made to depend
smoothly on $r$. Applying these diffeomorphisms to $A$ yields a smooth family
of submanifolds of $E$, parametrized by $r \gg 0$, which are Lagrangian with
respect to $\Omega^{(r)}$. One extends this arbitrarily to small values of $r$
to complete the definition of $\Lambda_{(c,[I])}$.

\begin{remarks}
(i) We will often write $c$, $\Lambda_c$ instead of $(c,[I])$,
$\Lambda_{(c,[I])}$, for the sake of brevity. This is fully justified only when
$n = 1$: in that case each component of the space of exact Lagrangian circles
on the surface $E_{c(1/2)}$ is simply-connected, so that $[I]$ contains no
information.

(ii) By definition each $\Lambda^{(r)}_c \subset E$, $r \in \Rgeq$, is
diffeomorphic to $A$. Since $A$ is glued together from two balls, it is a
homotopy $(n+1)$-sphere, but it may not be the standard one. To get standard
spheres from the construction, one has to impose the additional condition that
the composition
\[
\xymatrix{
 {S^n} \ar[r]^-{f_0} & {V_{c,0}} \ar[r] & {V_{c,1}} & \ar[l]_-{f_1} {S^n}
}
\]
can be isotoped within $\Diff(S^n)$ to an element of $O(n+1)$. Here $f_0,f_1$
are the framings of the vanishing cycles $V_{c,0},V_{c,1}$; and the middle
arrow is the diffeomorphism determined, up to isotopy, by $(I_t)$.
\end{remarks}

Assume now that no connected component of $S$ is closed, so that $\beta \in
\Omega^2(S)$ is exact. One then defines the Floer cohomology of two eventually
Lagrangian submanifolds $\Lambda,\Lambda' \subset \R^{\geq 0} \times E$ by
\begin{equation} \label{eq:eventually}
HF(\Lambda,\Lambda') = HF(\Lambda^{(r)},(\Lambda')^{(r)})
\end{equation}
where the right hand side is with respect to $\Omega^{(r)}$, for some $r \gg
0$. This makes sense because $\Omega^{(r)}$ is exact, and because there is a
natural class of compatible almost complex structures on $E$ such that
pseudo-holomorphic curves do not reach out to $\partial E$. Moreover, one can
show that the choice of $r$ does not affect \eqref{eq:eventually}.

\subsection{}
We will now explain the relation between the algebra and geometry introduced
above. Let $(E,\pi)$ and $\delta_{E/D}$ be an exact Morse fibration over $D$
with a relative Maslov map; the notation $M$ and $\delta_M$ is as usual. Let
$c$ be a matching path. One can find an admissible choice of paths
$(c_1,\dots,c_m)$ such that $c_i$ and $c_{i+1}$, for some $1 \leq i < m$, are
as in Figure \ref{fig:matching}, and such that no other $c_j$ intersects $c$.
In the corresponding distinguished basis $\Gamma = (L_1,\dots,L_m)$, $L_i$ and
$L_{i+1}$ are isotopic. Now choose gradings $\widetilde{\Gamma} =
(\tL_1,\dots,\tL_m)$ such that $\tL_i$ is isotopic to $\tL_{i+1}$, and let $\A
= \Lag(\widetilde{\Gamma})$ be the directed Fukaya category. Then
$(\tL_i,\tL_{i+1})$ is a matching pair in $D^b(\A)$. In fact one has
\begin{align*}
\Hom_{D^b(\A)}(\tL_i,\tL_{i+1}) & = H(hom_\A(\tL_i,\tL_{i+1})) \\ & =
HF^*(\tL_i,\tL_{i+1}) \iso H^*(L_i;\Z/2)
\end{align*}
which proves \eqref{eq:matching}, and the other property is similarly easy. The
informal idea is that the spherical object $C \in \Ob \, D^b(\A)$ associated to
this matching pair ``represents'', in some sense, the eventually Lagrangian
submanifold $\Lambda_c$. We will now make a concrete conjecture based on this
philosophy.
\includefigure{matching}{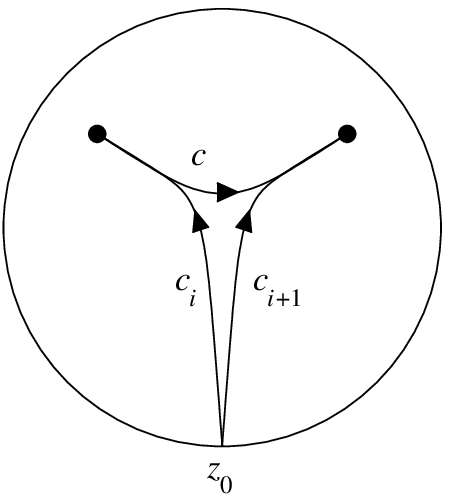}{hb}%

Suppose that $c,c'$ are two matching paths. For each of them one can find an
admissible choice of paths of the kind considered above, but in general these
choices cannot be the same for both paths. Therefore one gets two distinguished
bases $\Gamma,\Gamma'$ of vanishing cycles, and two directed Fukaya categories
$\A,\A'$ together with spherical objects $C \in \Ob \, D^b(\A)$ and $C' \in \Ob
\, D^b(\A')$, constructed from matching pairs of objects in $\A,\A'$ in the way
which we have explained. We know of course that $D^b(\A)$ and $D^b(\A')$ are
equivalent, but for the present purpose we need to make this slightly more
precise. There is a chain $\A = \A_0,\A_1,\dots,\A_r = \A'$ of directed
$A_\infty$-categories, of which each is related to the preceding one in one of
four possible ways: either $\A_{i+1}$ is quasi-isomorphic to $\A_i$; or it is
obtained from $\A_i$ by shifting the grading of the morphism spaces; or else
$\A_{i+1} = c\A$ or $r\A$. In the three last cases there are canonical
equivalences $D^b(\A_i) \iso D^b(\A_{i+1})$; in the first case we choose an
equivalence induced by an arbitrary quasi-isomorphism. Composing all of them
gives an equivalence $G: D^b(\A) \rightarrow D^b(\A')$.

\begin{conjecture} \label{th:induction}
Let $\Lambda_c,\Lambda_{c'}$ be the eventually Lagrangian submanifolds
associated to $c,c'$. Then
\begin{equation} \label{eq:induction}
HF(\Lambda_c,\Lambda_{c'}) \iso \bigoplus_{d \in \Z}
\Hom^d_{D^b(\A')}(G(C),C').
\end{equation}
\end{conjecture}

One can refine the statement by introducing suitable graded Floer cohomology
groups for $(\Lambda_c,\Lambda_{c'})$ but we prefer to skip this, since it does
not involve any really new ideas.

The two admissible choices of paths which lead to $\Gamma,\Gamma'$ are not
unique, and neither is the sequence of Hurwitz moves joining them; moreover, in
the definition of $G$ we admitted arbitrarily chosen quasi-isomorphisms. It is
implicitly part of the conjecture that the right hand side of
\eqref{eq:induction} is independent of all these choices. The most
doubtful-looking point is the bit about quasi-isomorphisms. Imagine for
instance that $\A_i$ and $\A_{i+1}$ are directed Fukaya categories arising from
the same Lagrangian configuration, but with different choices of almost complex
structures. As mentioned in \cite[Section 6]{seidel00}, there is an explicit
construction of quasi-isomorphisms between $\A_i$ and $\A_{i+1}$, using
pseudo-holomorphic curves for a one-parameter family of choices of almost
complex structures. One might think that for \eqref{eq:induction} to be
correct, it would be necessary to choose the step between $D^b(\A_i)$ and
$D^b(\A_{i+1})$ in the construction of $G$ to be such an analytically defined
quasi-isomorphism. The answer to this objection is Lemma \ref{th:invariant},
which says that choices of quasi-isomorphisms can never affect how objects like
$C$ are mapped.

To see what Conjecture \ref{th:induction} does, consider the case $n = 1$. Then
$M$ is a surface, so that $\A$ and $\A'$ can be determined combinatorially, and
so can $G: \A \rightarrow \A'$ (this is only true thanks to the freedom of
choosing quasi-isomorphisms, which is why we were insisting on it). On the
other hand, the left hand side of \eqref{eq:induction} is the Floer cohomology
of Lagrangian two-spheres in the symplectic four-manifold $(E,\Omega^{(r)})$,
which is difficult or even impossible to compute directly.

\begin{remark} \label{th:zero-spheres}
With some modifications, the discussion can be adapted to the trivial dimension
$n = 0$. Here we only want to point out one example, namely the double branched
cover $\pi: E \rightarrow D$ from Section \ref{sec:dimension-zero}c. Because
the smooth fibre consists of two points, any path is a matching path, so that
there is a multitude of one-dimensional spherical objects. This explains the
geometric meaning of the objects $C_1,\dots,C_{2g} \in \Ob \, D^b(\A_g)$ which
we wrote down earlier. By pushing the argument further, one can show that the
zero-dimensional version of Conjecture \ref{th:induction} implies the algebraic
formulae for geometric intersection numbers on $E$ which were derived, by a
quite different argument, in \cite{khovanov-seidel98}.
\end{remark}

%
%
\providecommand{\bysame}{\leavevmode\hbox to3em{\hrulefill}\thinspace}


\begin{thebibliography}{VC}
\bibitem[VC]{seidel00}
P.~Seidel, \emph{Vanishing cycles and mutation}, Proceedings of the 3rd
European Congress of Mathematics (Barcelona, 2000), Birkh{\"a}user, to appear.

\bibitem{acampo00}
N.~A'Campo, \emph{A combinatorial property of generic immersions of curves},
  Preprint.

\bibitem{acampo73b}
\bysame, \emph{Le nombre de lefschetz d'une monodromie}, Indag. Math.
  \textbf{35} (1973), 113--118.

\bibitem{acampo75b}
\bysame, \emph{Le groupe de monodromie du d{\'e}ploiement des singularit{\'e}s
  isol{\'e}es de courbes planes {I}}, Math. Ann. \textbf{213} (1975), 1--32.

\bibitem{acampo98}
\bysame, \emph{Real deformations and complex topology of plane curve
  singularities}, Annales Fac. Sci. Toulouse \textbf{8} (1999), 5--23.

\bibitem{arnold-gusein-zade-varchenko}
V.~I. Arnold, S.~M. Gusein-Zade, and A.~N. Varchenko, \emph{Singularities of
  differentiable maps}, Birkh{\"a}user, 1988.

\bibitem{beilinson78}
A.~Beilinson, \emph{Coherent sheaves on {$\mathbb{P}^n$} and problems of
  linear algebra}, Funct. Anal. Appl. \textbf{12} (1978), 214--216.

\bibitem{bernstein-gelfand-ponomarev72}
J.~Bernstein, I.~Gelfand, and V.~Ponomarev, \emph{Coxeter functors and
  {G}abriel's theorem}, Russian Math. Surveys \textbf{28} (1973), 17--32.

\bibitem{bernstein-lunts94}
J.~Bernstein and V.~Lunts, \emph{Equivariant sheaves and functors}, Lecture
  Notes in Mathematics, vol. 1578, Springer, 1994.

\bibitem{bondal-polishchuk94}
A.~Bondal and A.~Polishchuk, \emph{Homological properties of associative
  algebras: the method of helices}, Russian Math. Izvestiya \textbf{42} (1994),
  219--260.

\bibitem{cohen-jones-segal95}
R.~L. Cohen, J.~D.~S. Jones, and G.~B. Segal, \emph{{F}loer's
  infinite-dimensional {M}orse theory and homotopy theory}, The {F}loer
  memorial volume (H.~Hofer, C.~Taubes, A.~Weinstein, and E.~Zehnder, eds.),
  Progress in Mathematics, vol. 133, Birkh{\"a}user, 1995, pp.~297--325.

\bibitem{fukaya97}
K.~Fukaya, \emph{Floer homology for three-manifolds with boundary {I}},
  Preprint, 1997.

\bibitem{givental94}
A.~Givental, \emph{Homological geometry and mirror symmetry}, Proceedings of
  the International Congress of Mathematics, Z{\"u}rich, vol.~1,
  Birkh{\"a}user, 1994, pp.~472--480.

\bibitem{gusein-zade74}
S.~M. Gusein-Zade, \emph{Dynkin diagrams for singularities of functions of two
  variables}, Functional Anal. Appl. \textbf{8} (1974), 295--300.

\bibitem{hori-iqbal-vafa00}
K.~Hori, A.~Iqbal, and C.~Vafa, \emph{{$D$}-branes and mirror symmetry},
  Preprint hep-th/0005247.

\bibitem{kapranov90}
M.~M. Kapranov, \emph{{M}utations and {S}erre functors on constructive
  bundles}, Functional Anal. Appl. \textbf{24} (1990), 155--156.

\bibitem{khovanov-seidel98}
M.~Khovanov and P.~Seidel, \emph{Quivers, {F}loer cohomology, and braid group
  actions}, Preprint Math.QA/0006056.

\bibitem{koenig-zimmermann98}
S.~Koenig and A.~Zimmermann (eds.), \emph{Derived equivalences for group
  rings}, Lecture notes in Math., vol. 1685, Springer, 1998.

\bibitem{le74}
D.~T. L{\^e}, \emph{La monodromie n'a pas de points fixes}, J. Fac. Sci. Univ.
  Tokyo Sect. IA \textbf{22} (1975), 409--427.

\bibitem{le-ono97}
H.~V. L{\^e} and K.~Ono, \emph{Cup-length estimates for symplectic fixed
points},
  Contact and symplectic geometry (Cambridge, 1994), Cambridge Univ. Press,
  1996, pp.~268--295.

\bibitem{penkava-schwarz95}
M.~Penkava and A.~Schwarz, \emph{{$A_\infty$}-algebras and the cohomology of
  moduli spaces}, Amer. Math. Soc. Translations, vol. 169, pp.~91--108,
  American Mathematical Society, 1995.

\bibitem{schwarz96}
M.~Schwarz, \emph{A quantum cup-length estimate for symplectic fixed points},
  Invent. Math. \textbf{133} (1998), 353--397.

\bibitem{seidel96b}
P.~Seidel, \emph{The symplectic {F}loer homology of a {D}ehn twist}, Math.
  Research Lett. \textbf{3} (1996), 829--834.

\bibitem{seidel-thomas99}
P.~Seidel and R.~Thomas, \emph{Braid group actions on derived categories of
  coherent sheaves}, Preprint math.AG/0001043. To appear in {\em Duke Math.\
  J}.

\bibitem{wajnryb83}
B.~Wajnryb, \emph{A simple presentation for the mapping class group of an
  orientable surface}, Israel J. Math. \textbf{45} (1983), 157--174.
\end{thebibliography}
\end{document}